\documentclass[a4paper, 11pt, twoside]{amsart}
\author[Antonio Rapagnetta]{Antonio Rapagnetta}
\address{Dipartimento di Matematica universit\`a di Roma `Tor
Vergata'.\newline Via della Ricerca Scientifica 00133 Roma} \email{rapagnet@mat.uniroma2.it}
\title[Topological invariants of $\widetilde{\mathcal{M}}$]
{Topological invariants of O'Grady's six dimensional irreducible symplectic variety}

\usepackage{amsmath}
\usepackage{amsfonts}
\usepackage{amscd}
\usepackage{amsthm}
\usepackage{graphics}
\usepackage{verbatim}

\newtheorem{thm}{Theorem}[subsection]
\newtheorem{prop}[thm]{Proposition}
\newtheorem{lem}[thm]{Lemma}
\newtheorem{cor}[thm]{Corollary}
\theoremstyle{definition}
\newtheorem{defn}[thm]{Definition}
\newtheorem{claim}[thm]{Claim}
\newtheorem{notation}[thm]{Notation}

\newtheorem{rem}[thm]{Remark}

\begin{document}
\begin{abstract}
We study O'Grady examples of irreducible symplectic varieties: we establish that both of them can be deformed into lagrangian fibrations. We analyse in detail the topology of the six dimensional example: in particular we compute its Euler characteristic and determine its Beauville form.
\end{abstract}
\maketitle
\section*{Introduction}
 Irreducible symplectic varieties are defined as compact
K\"ahler varieties having trivial fundamental group and endowed
with a unique global holomorphic 2-form which is non degenerate on
each point.\\
By the Bogomolov decomposition \cite{Bo74}, irreducible
symplectic varieties play  (together with  Calabi Yau manifolds
and Complex tori) a central role in the classification of compact
K\"ahler
manifolds with torsion $c_{1}$.\\
Very few examples  of irreducible symplectic varieties are
available in  literature.\\For any positive integer $n$ Beauville
exhibited 2 examples of dimension $2n$ (\cite{Be83}): the Hilbert
scheme $Hilb^{n}(X)$ parametrizing 0 dimensional subschemes of
length $n$ on a K3 surface $X$, and the Kummer generalized variety
$K^{n}(T)$ of a 2-dimensional torus $T$, namely the locus in
$Hilb^{n+1}(T)$
parametrizing subschemes having associated cycle summing up to zero.\\
Besides the Beauville examples, there are only two known examples
of irreducible symplectic varieties up to deformation equivalence:
they have been exhibited  by O'Grady in \cite{OG99} and \cite{OG03} and their dimensions are
respectively ten and six.\\
While Hilbert schemes of points and Kummer generalized varieties
have been deeply studied, very little is known about the two
O'Grady examples.\\
In this paper, after having proved that both the O'Grady examples
can be deformed into Lagrangian fibrations (see Corollary
\ref{lagrang}) we study the topology of the one,
$\widetilde{\mathcal{M}}$, having dimension six. In section 2 we
establish that its Euler characteristic is 1920 (see Theorem
\ref{caratteristica m tilde}). Finally, in section 3 we determine
the Beauville form and the Fujiki
constant of $\widetilde{\mathcal{M}}$ (see Theorem \ref{mainthm}).\\
 It is remarkable that the Fujiki constant of
$\widetilde{\mathcal{M}}$ is $60$, as in the case of  the
generalized Kummer variety of the same dimension: this is
the first known case where two non diffeomorphic irreducible
symplectic varieties have the same Fujiki constant.\\
{\bf Acknowledgements.} This paper contains results taken from my
PhD thesis. It is a pleasure to thank my advisor Kieran O'Grady
for having followed me during the elaboration of the thesis, and
more generally for all the mathematics that he has taught me.
\section{O'Grady's desingularization}

In this section  we recall the construction of O'Grady's examples
of irreducible symplectic varieties and we give a slight
generalization of O'Grady's symplectic desingularization (see
propositions \ref{desingk3} and \ref{desing abel}). This
desingularization enables us to prove in Corollary \ref{lagrang}
that all known examples of irreducible symplectic varieties are
deformation equivalent to lagrangian fibrations. In the rest of
the section we fix the notation that we will use in the next 2
sections where we study the six dimensional O'Grady example
$\widetilde{\mathcal{M}}$.

\subsection{Analysis of O'Grady's construction}
 We need to recall two general
definitions due to Mukai. In these definitions $X$ is a projective
symplectic surface, namely $X$ is a $K3$ or an abelian surface.
\begin{defn}\label{def reticolo mukai}Consider
the involution on the even cohomology of $X$ given by
\[\begin{array}{ccc}
H^{ev}(X,\mathbb{Z})& \longrightarrow &
H^{ev}(X,\mathbb{Z})\\
\alpha=\alpha_{0}+\alpha_{2}+\alpha_{4} & \mapsto &
\overline{\alpha}=\alpha_{0}-\alpha_{2}+\alpha_{4}\\
\end{array}
\]
with $\alpha_{2i}\in H^{2i}(X,\mathbb{Z})$. Let $\alpha ,\beta\in
H^{ev}(X,\mathbb{Z})$: the Mukai's pairing $<\cdot,\cdot>$ is
given by:
\[<\alpha,\beta>:=-\int\overline{\alpha}\cup\beta
.\] Finally the Mukai lattice is
$(H^{ev}(X,\mathbb{Z}),<\cdot,\cdot>)$.
\end{defn}
\begin{defn}\label{definizione vettore di mukai}
Let $F$ be a coherent sheaf on $X$, the Mukai vector of $F$ is
\[m(F):=ch(F)\cup\sqrt{Td(X)}\in H^{ev}(X,\mathbb{Z}).\]
\end{defn}
\begin{notation}\label{notazione moduli}
Letting $\eta\in H^{4}(X,\mathbb{Z})$ be the fundamental class,
given
\[v=(v_{0},v_{2},v_{4}\eta)\in H^{0}(X,\mathbb{Z})\oplus
\mathbf{Pic}(X)\oplus H^{4}(X,\mathbb{Z}),\] we will denote by
$\mathbf{M}_{v}$ the Simpson moduli space of semistable sheaves on
$X$ having  as Mukai's vector the class of $v$ in $
H^{0}(X,\mathbb{Z})\oplus H^{2}(X,\mathbb{Z})\oplus
H^{4}(X,\mathbb{Z})$.\\
Using the identification $H^{4}(X,\mathbb{Z})\equiv\mathbb{Z}$ we
will usually denote $v=(v_{0},v_{2},v_{4}\eta)$ by
$(v_{0},v_{2},v_{4})$. Moreover, to simplify the notation, we will
always replace $v_{2}$ by a divisor in the linear equivalence
class defined by $v_{2}$.
\end{notation}
 In this section we simply fix the notation and briefly recall
the construction of O'Grady's examples of irreducible symplectic
varieties.\\
 In all this paper ${\mathcal{J}}$ will be the
Jacobian of a genus 2 curve $C^{0}$ such that
$\mathbf{NS}(\mathcal{J})=\mathbb{Z}c_{1}(\Theta)$ where $\Theta $
is
 a symmetric theta divisor .\\
According to Notation \ref{notazione moduli}, given
\[v=(v_{0},v_{2},v_{4})\in H^{0}(\mathcal{J},\mathbb{Z})\oplus
\mathbf{Pic}(\mathcal{J})\oplus H^{4}(\mathcal{J},\mathbb{Z}),\]
we will denote by $\mathbf{M}_{v}$ the Simpson moduli space of
semistable sheaves on $\mathcal{J}$ having  as Mukai's vector the
class of $v$ in $ H^{ev}(\mathcal{J},\mathbb{Z})$.\\
The moduli space $\mathbf{M}_{v}$ is endowed with a regular
morphism defined by
\begin{equation}\label{definizione a}
\begin{matrix}
a_{v}: &\mathbf{M}_{v}&\longrightarrow
&\mathcal{\widehat{J}}\times\mathcal{J}\\
&[F]&\longmapsto &(Det(F),\sum
c_{2}(F))\\
\end{matrix}\end{equation}
where $[F]$ is the S-equivalence class of $F$, $Det(F)$ is the
determinant bundle of $F$ and, if the formal sum $\sum n_{i}p_{i}$
is a representative of the second Chern class of $F$ in the Chow
ring of $\mathcal{J}$, then $\sum
c_{2}(F):=\sum n_{i}p_{i}\in\mathcal{J}$ .\\
This enables us to give the following definition
\begin{defn}
\[\mathbf{M}_{v}^{0}:=a_{v}^{-1}(v_{2},0).\]
\end{defn}
Obviously the points of $\mathbf{M}_{v}^{0}$ parametrize
S-equivalence classes  of semistable sheaves having determinant
bundle linearly equivalent to $v_{2}$ and second Chern class
summing up to $0$.
\\In  the paper \cite{OG03} O'Grady constructed his
second example of irreducible symplectic variety starting from the
moduli space $\mathbf{M}_{(2,0,-2)}$.\\The singular locus $\Sigma$
of $\mathbf{M}_{(2,0,-2)}$ parametrizes polystable sheaves of the
form \[I_{p}\otimes L_{1}\oplus I_{q}\otimes L_{2}\] where $p$ and
$q$ are points on $\mathcal{J}$, $I_{p}$ and $I_{q}$ are the
corresponding sheaves of ideals, and $L_{i}$ are line
bundles with homologically trivial first Chern class.\\
The singular locus $\Omega$ of $\Sigma$ is precisely  the
subscheme parametrizing   sheaves
\[I_{p}\otimes L\oplus
I_{p}\otimes L.\] The first step in constructing the O'Grady
6-dimensional  irreducible symplectic variety is the following:
blow up $\mathbf{M}_{(2,0,-2)}$ along $\Omega$ and then blow up
$Bl_{\Omega}\mathbf{M}_{(2,0,-2)}$ along the strict transform of
$\Sigma$ and finally contract the inverse image of $\Omega$ via
the two blow ups. This  produces  a ten dimensional smooth variety
$\widetilde{\mathbf{M}}_{(2,0,-2)}$. Moreover the obvious
birational map from $\widetilde{\mathbf{M}}_{(2,0,-2)}$ to
$\mathbf{M}_{(2,0,-2)}$, extends to a regular morphism
\[\widetilde{\pi}: \widetilde{\mathbf{M}}_{(2,0,-2)}\longrightarrow\mathbf{M}_{(2,0,-2)}.\]
Let  $a_{(2,0,-2)}$ be defined by (\ref{definizione a}), and
define
\begin{equation}
\widetilde{a}_{(2,0,-2)}:=a_{(2,0,-2)}\circ\widetilde{\pi}:
\widetilde{\mathbf{M}}_{(2,0,-2)}\longrightarrow
\mathcal{\widehat{J}}\times\mathcal{J}.
\end{equation} The new  six dimensional irreducible symplectic variety of O'Grady is then
given by:\begin{defn}\label{definizione mtilde}
\[\widetilde{\mathcal{M}}:=\widetilde{a}_{(2,0,-2)}^{-1}(0,0).
\]\end{defn}
\begin{rem}\label{costruzione og rem}
$\widetilde{\mathcal{M}}$ can be equivalently constructed starting
from the locus $\mathbf{M}_{(2,0,-2)}^{0}=a_{(2,0,-2)}^{-1}(0,0)$
blowing up $\mathbf{M}_{(2,0,-2)}^{0}$ along $\Omega\cap
\mathbf{M}_{(2,0,-2)}^{0}$, blowing up here the strict transform
of $\Sigma\cap \mathbf{M}_{(2,0,-2)}^{0}$ and contracting the
inverse image via the composition of the two blow ups of
$\Omega\cap \mathbf{M}_{(2,0,-2)}^{0}$.
\end{rem}
\begin{rem}There are Mukai vectors different from $(2,0,-2)$ for
which  O'Grady's construction in \cite{OG99}, sketched above,
works and produces a smooth algebraic variety with a holomorphic
symplectic two form. They probably do not give new deformation
classes of irreducible symplectic varieties, but at least two of
them will be useful in this paper to understand the geometry of
O'Grady's examples, so in the next two propositions we state which
Mukai vectors admit such a weak generalization of O'Grady's result
and in the successive corollary we  single out the cases that we
will effectively use later.
\end{rem}
\begin{prop}\label{desingk3}
Let $X$ be a K3 surface or an abelian surface such that
$\mathbf{NS}(X)=\mathbb{Z}H$. Let $v\in H^{0}(X,\mathbb{Z})\oplus
\mathbf{Pic}(X)\oplus H^{4}(X,\mathbb{Z})$ and
let $\overline{v}\in H^{ev}(X,\mathbb{Z})$ be the class of $v$.\\
Suppose that:
\begin{enumerate}\item{$2|\overline{v}$, but $\frac{\overline{v}}{2}$ is
primitive,} \item{$<\overline{v},\overline{v}>=8$,} \item{The
moduli space $\mathbf{M}_{\frac{v}{2}}$ is fine and non
empty.}\end{enumerate} Then $\mathbf{M}_{v}$ is reduced and there
exists a symplectic desingularization
\[\widetilde{\pi}_{v}:\widetilde{\mathbf{M}}_{v}\rightarrow \mathbf{M}_{v}.\]
It can be obtained exactly repeating O'Grady's construction in
\cite{OG99}.
\end{prop}
\begin{rem}
 A slight modification of the following proof shows that
Proposition \ref{desingk3} still holds if we only require in 3)
$\mathbf{M}_{\frac{v}{2}}\neq \emptyset$. \\
This condition is always verified if the first non zero
coefficient of $\frac{v}{2}$ is positive (see Theorem 0.1 of
\cite{Yo} and Theorem 0.1 of \cite{Yo01}).
\end{rem}
\begin{proof}
The proof is exactly the same after replacing the sheaves of the
form $I_{p}\otimes L$  by the ones whose Mukai vector is
$\frac{1}{2}\overline{v}$. Since $\mathbf{NS}(X)=\mathbb{Z}H$, the
Hilbert polynomial of a sheaf determines its Mukai vector, hence,
by primitivity, $\mathbf{M}_{\frac{v}{2}}$ is smooth and a
strictly semistable sheaf $F$ such that  $[F]\in\mathbf{M}_{v}$
fits in an exact sequence
\begin{equation}\label{*}0\rightarrow G_{1}\rightarrow F\rightarrow
G_{2}\rightarrow 0\end{equation} with
$[G_{i}]\in\mathbf{M}_{\frac{v}{2}}$.\\
For $\overline{v}$ divisible only by $2$, the classification of
the structures of semistable sheaves with their automorphism
groups modulo scalars, is easily seen to be the following:
\begin{enumerate}
\item{$\frac{Aut(F)}{\mathbb{C}^{*}}=PGL(2)$ if $G_{1}=G_{2}$ and
the extension (\ref{*}) is trivial,}
\item{$\frac{Aut(F)}{\mathbb{C}^{*}}=(\mathbb{C},+)$ if
$G_{1}=G_{2}$ and the extension (\ref{*}) is non
trivial,}\item{$\frac{Aut(F)}{\mathbb{C}^{*}}=\mathbb{C}^{*}$ if
$G_{1}\neq G_{2}$ and the extension (\ref{*}) is  trivial,}
\item{$\frac{Aut(F)}{\mathbb{C}^{*}}=id$ if $G_{1}\neq G_{2}$ and
the extension (\ref{*}) is non trivial.}
\end{enumerate}
In each of these items $[G_{i}]\in\mathbf{M}_{\frac{v}{2}}$: this
generalizes
 Corollary (1.1.8) of \cite{OG99}.\\
In section (1)  of \cite{OG99} for any even $c\ge 4$ a
desingularization $\widehat{\mathbf{M}}_{(2,0,2-c)}$ of
$\mathbf{M}_{(2,0,2-c)}$ is constructed. Any statement proved in
section (1)  of \cite{OG99} has an obvious analogous if we
replace the Mukai vector $(2,0,2-c)$ with a vector $\overline{v}$
satisfying our hypothesis.  Furthermore, using the given
classification of semistable sheaves, we can repeat   exactly the
same proofs  if we assume the followings:
\begin{itemize}
\item{$\mathbf{M}_{\frac{v}{2}}$ is endowed with a tautological
family; this property of $\mathbf{M}_{(1,0,1-\frac{c}{2})}$ is
used in the preparation and in the proof of  Proposition
(1.7.10),}
\item{for $[G_{1}]\neq [G_{2}]$ in $\mathbf{M}_{\frac{v}{2}}$,
there are non trivial extensions of $G_{1}$ by $G_{2}$; this
property for $[G_{i}]\in \mathbf{M}_{(1,0,1-\frac{c}{2})}$ is used
in the proof of Claim (1.4.8),}\item{$dim(Ext^{1}(G_{i},G_{i}))\ge
4$ for any $[G_{i}]\in\mathbf{M}_{\frac{v}{2}}$; this property for
$G_{i}\in \mathbf{M}_{(1,0,1-\frac{c}{2})}$ is used in the proof
of Lemma (1.5.6).}
\end{itemize}
Since \[dim(Ext^{1}(G_{1},G_{2}))=\{\begin{array}{cc} \frac{<\overline{v},\overline{v}>}{4} & if\; G_{1}\neq G_{2}\\
\frac{<\overline{v},\overline{v}>}{4}+2 & if \; G_{1}=G_{2}
\end{array}\]
we get that under our hypothesis there exists a desingularization
$\widehat{\mathbf{M}}_{v}$ of $\mathbf{M}_{v}$ obtained repeating
formally O'Grady's desingularization of $\mathbf{M}_{(2,0,2-c)}$.
In  section (2) O'Grady fixes  $c=4$ in order to obtain a
contraction of $\widehat{\mathbf{M}}_{(2,0,2-c)}$ to a symplectic
variety $\widetilde{\mathbf{M}}_{(2,0,2-c)}$. For $c=4$ the
desingularization procedure is more simple: the simplification
depends only on \[dim(Ext^{1}(G,G))=4 \;\forall G\in
\mathbf{M}_{(1,0,1-\frac{c}{2})}\](see Formula (1.8.2) and the
successive paragraph). In a completely analogous way, if
\begin{equation}\label{ipotesi} dim(Ext^{1}(G,G))=4 \;\forall G\in
\mathbf{M}_{\frac{v}{2}}\end{equation} the desingularization of
$\mathbf{M}_{v}$ simply consists in blowing up $\mathbf{M}_{v}$
along the closed subvariety $\Omega_{v}$ whose points represent to
sheaves of the form $G^{2}$ and then blowing up the strict
transform of the closed subvariety  $\Sigma_{v}$ whose points
represent all the semistable sheaves.\\
Since hypothesis 3) implies (\ref{ipotesi}) and since any result
in section 2) of \cite{OG99} is a consequence of the analysis in
section 1) we also get that those  results  still hold if we
replace $\mathbf{M}_{(2,0,-2)}$ by $\mathbf{M}_{v}$ with $v$
stasfying our hypotheses. In particular there exists a symplectic
birational model $\widetilde{\mathcal{M}}_{v}$ of $\mathbf{M}_{v}$
and it is endowed with a regular map
\[\widetilde{\pi}_{v}:\widetilde{\mathbf{M}}_{v}\rightarrow
\mathbf{M}_{v}\] being an isomorphism on the smooth locus of
$\mathbf{M}_{v}$.\\
It remains to prove that $\mathbf{M}_{v}$ is reduced. From section
(1)  of \cite{OG99} we can directly deduce that
$\mathbf{M}_{(2,0,2-c)}$ is reduced. In fact
$\mathbf{M}_{(2,0,2-c)}$ is the $PGL(N)$-quotient of the scheme
$\mathcal{Q}_{c}$, so it will be enough to show that the
semistable locus of $\mathcal{Q}_{c}$ is reduced. Since any point
representing a semistable sheaf can be moved by $PGL(N)$ to any
neighborhood of any point representing the polystable sheaf
associated to it, it will be enough to prove that any point in
$\mathcal{Q}_{c}$ representing a polystable sheaf has a reduced
neighborhood. Neighborhoods in $\mathcal{Q}_{c}$  of points
representing  polystable sheaves are described in section (1) of
\cite{OG99} . If the sheaf $F$ represented by the point $y$ is
polystable, the $PGL(n)-$orbit of $y$ is closed, and applying
Luna's etale  slice theorem (see pages 54-55 of \cite{OG99}) we
get that a neighborhood of $y$ has an etale covering from
$PGL(n)\times_{st(y)}\mathcal{V}$, where $st(y)$ is the stabilizer
of $y$ and $\mathcal{V}$ is the etale slice. In the cases
occurring in the study of $\mathcal{Q}_{c}$, we always get that
$\mathcal{V}$ is reduced, locally irreducible near $y$. In fact,
if $\mathcal{V}$ is not smooth near $y$ (namely if $F$ is not
stable), the normal cone to $\mathcal{V}$ at $y$ turns out to be
either a the affine cone over a reduced irreducible  quadric
hypersurface if $F=G_{1}\oplus G_{2}$ with
$G_{1}\neq G_{2}$ (see Claim (1.4.8) and Proposition (1.4.10)) or
the affine cone over a reduced irreducible complete intersection
of three quadrics (see Lemma (1.5.6) and Proposition (1.5.10)). It
follows that $PGL(n)\times_{st(y)}\mathcal{V}$ is always reduced
near $(1,y)$: hence its etale image
in $\mathcal{Q}_{c}$ is reduced too.\\
Since as we said earlier any statement of  section 1 of \cite{OG99}
 still holds, with the same proof, after replacing  $(2,0,2-c)$
with a $v$ satisfying our hypotheses, then we get that for any
such a $v$ the moduli space $\mathbf{M}_{v}$ is reduced.
\end{proof}
\begin{rem}\label{fibre ogradi}
In the proof of the previous proposition we showed that all that
is proved in sections 1) and 2) of \cite{OG03} holds under the
hypothesis of the proposition. In particular we can describe the
fiber of O'Grady's desingularization as follows. If $p\in
\mathbf{M}_{v}$ corresponds to a sheaf $F$ such that
$\frac{Aut(F)}{\mathbb{C}^{*}}=PGL(2)$ then
$\widetilde{\pi}_{v}^{-1}(p)$ is a smooth 3-dimensional quadric as
in formula (2.2.9) on page 88 of \cite{OG99}. If $p\in
\mathbf{M}_{v}$ corresponds to a sheaf $F$ such that
$\frac{Aut(F)}{\mathbb{C}^{*}}=\mathbb{C}^{*}$ then
$\widetilde{\pi}_{v}^{-1}(p)$ is a $\mathbb{P}^{1}$ as shown in
formula (2.2.4) on page 87 of \cite{OG99}. This will be used in
the computation of the Euler characteristic of
$\widetilde{\mathcal{M}}$.\end{rem}
\begin{prop}\label{desing abel} Let $X=\mathcal{J}$,
let $v$ satisfy the hypothesis of Proposition \ref{desingk3}. Let
$\widetilde{\pi}:
\widetilde{\mathbf{M}}_{v}\longrightarrow\mathbf{M}_{v}$ be the
desingularization map obtained in the same proposition and let
$\widetilde{\mathbf{M}}_{v}^{0}$ be as defined in \ref{definizione
mtilde}, then $\widetilde{\mathbf{M}}_{v}^{0}$ is a smooth
algebraic variety endowed with a symplectic holomorphic two form
obtained restricting the one on $\widetilde{\mathbf{M}}_{v}$.
\end{prop}
\begin{proof}We have only to prove that the symplectic form
of $\widetilde{\mathbf{M}}_{v}$ restricts to a symplectic form on
$\widetilde{\mathbf{M}}_{v}^{0}$: this proof can be copied from
Proposition 1.1 of \cite{De99} and Proposition 2.3.3 of
\cite{OG03}.
\end{proof}
\begin{rem}Notice that, when Proposition \ref{desingk3} and
Proposition \ref{desing abel} actually work, namely when the
moduli spaces involved are not empty, they always produce
respectively pure 10 dimensional and pure 6 dimensional symplectic
varieties.
\end{rem}
\begin{cor}\label{vettori espliciti}\begin{enumerate} \item{Set $X:=\mathcal{J}$, the
procedure described in the previous propositions produces a
symplectic desingularization $\widetilde{\mathbf{M}}_{(0,
2\Theta,-2)}^{0}$ of $\mathbf{M}_{(0, 2\Theta,-2)}^{0}$.}
\item{Let $X$ be a K3 surface obtained as a double covering of the
projective plane ramified along a smooth sextic. Suppose  that
$\mathbf{Pic}(X)$ is generated by $H$, the pull-back of a line.
Then there exists a symplectic desingularization
$\widetilde{\mathbf{M}}_{(0,2H,2)}$ of $\mathbf{M}_{(0,2H,2)}$.}
\end{enumerate}
\end{cor}
\begin{proof}
Conditions (1) and (2) of Proposition \ref{desingk3} are obviously
satisfied, the condition (3) is satisfied using the criterion (see
appendix of \cite{Mu84}) asserting that a tautological family for
the stable locus of a moduli space $\mathbf{M}_{v}$ exists when
$G.C.D(v_{0},\{v_{2}\cdot c_{1}(L)\}_{L\in \mathbf{Pic}(X)},
\chi)=1$ ($\chi$ is the Euler characteristic of any sheaf of
$\mathbf{M}_{v}$), indeed in the cases we are considering $\chi
=-1$ or $\chi =1$.
\end{proof}
\begin{rem} \label{non solo symplettica}In Proposition \ref{v02-2}
 we will prove that $\widetilde{\mathbf{M}}_{(0,2\Theta,-2)}^{0}$
 is birational to $\widetilde{\mathcal{M}}$, on
the other hand in Proposition (4.1.5) of \cite{OG99} it is proved
(with a different notation) that the 10 dimensional O'Grady's
example is birational to an irreducible component $N$ of
$\mathbf{M}_{(0,2H,2)}$ (actually it can be proved that
$\mathbf{M}_{(0,2H,2)}$ is irreducible). Since simple connectivity
and $dim(H^{2,0})$ are birational invariants for smooth varieties,
this  implies that $\widetilde{\mathbf{M}}_{(0,2\Theta,-2)}^{0}$
and  the irreducible  component $\widetilde{N}$ (over $N$) of
$\widetilde{\mathbf{M}}_{(0,2H,2)}$ are symplectic irreducible.
\end{rem}
It is conjectured  (see \cite{Sa03}) that any Irreducible
symplectic variety is deformation equivalent to a Lagrangian
fibration (i.e. a 2n dimensional symplectic variety endowed with a
proper morphism to an n-dimensional variety, such that the
restriction of the symplectic form to the general fiber is zero):
in the following corollary we verify this conjecture on the
known examples.
\begin{cor}\label{lagrang}
All the known irreducible symplectic varieties are deformation
equivalent to  lagrangian fibrations.
\end{cor}
\begin{proof}
It is well known for Hilbert schemes of points on a K3 surfaces
(see \cite{Be99}) and for generalized Kummer varieties (see
\cite{De99}). For the O'Grady examples, by the last remark and
since birational irreducible symplectic varieties are deformation
equivalent (see \cite{Hu97}), it is enough to prove that
$\widetilde{\mathbf{M}}_{(0,2\Theta,-2)}^{0}$ and $\widetilde{N}$
are lagrangian fibrations. Both of these spaces are
desingularizations of moduli spaces parametrizing sheaves on
surfaces  supported in codimension 1 : the functor associating to
a sheaf the subscheme defined by its Fitting ideal (see
\ref{fitting astra})always defines, on any such moduli space, a
regular map to the suitable Hilbert scheme. In the case of
$N\subset\mathbf{M}_{(0, 2H, 2)}^{0}$ the Hilbert scheme is
identified with the 5 dimensional linear system $|2H|$ and the
morphism is easily seen (see 4.2 of \cite{OG99}) to
be surjective.\\
In the case of $\mathbf{M}_{(0, 2\Theta,-2)}^{0}$ we will see in
section 3.2.4 that this regular morphism (we will denote this map
by $\Phi$) is surjective on a closed subvariety of the Hilbert
scheme identified with the 3 dimensional linear
system $|2\Theta|$.\\
Therefore both the irreducible symplectic varieties
$\widetilde{\mathbf{M}}_{(0,2\Theta,-2)}^{0}$ and $\widetilde{N}$
have a surjective morphism to a projective space of dimension
equal to half their dimension: by Theorem 1 of \cite{Ma01} these
morphisms are lagrangian fibrations.
\end{proof}
\subsection{Notation}Now we fix the notations that we will follow for
the rest of the paper.
\begin{notation}\label{grandi notazioni}
Given  a coherent sheaf $F$ on our abelian surface $\mathcal{J}$,
we denote by $[F]$
its S-equivalence class or equivalently the associated point in the moduli space.\\
Given $v\in  H^{0}(\mathcal{J},\mathbb{Z})\oplus
\mathbf{Pic}(\mathcal{J})\oplus H^{4}(\mathcal{J},\mathbb{Z})$
satisfying the hypotheses of  Proposition \ref{desing abel},
denoting by $\overline{v}$ its class in $H^{ev}(\mathcal{J},\mathbb{Z})$, we set:\\
$\Sigma_{v}:=\{[F_{1}\oplus F_{2}]\in \mathbf{M}_{v}:
m(F_{1})=ch(F_{1})= m(F_{2})= ch(F_{2})=\frac{1}{2}\overline{v}\}$
the singular locus
of $\mathbf{M}_{v}$,\\
$\Omega_{v}:=\{[F\oplus F]\in \mathbf{M}_{v}:
m(F)=ch(F)=\frac{1}{2}\overline{v} \}$ the singular locus of $\Sigma_{v}$,\\
$\widetilde{\pi}_{v}:\widetilde{\mathbf{M}}_{v}\longrightarrow
\mathbf{M}_{v}$  the symplectic desingularization
map,\\
$\widetilde{\Sigma}_{v}:=\widetilde{\pi}_{v}^{-1}\Sigma_{v}$,\\
$\widetilde{\Omega}_{v}:=\widetilde{\pi}_{v}^{-1}\Omega_{v}$,\\
$a_{v}:=Det
\times\Sigma c_{2}:\mathbf{M}_{v}\longrightarrow
\mathbf{Pic}(\mathcal{J})\times\mathcal{J}$,\\
$\widetilde{a}_{v}:=a_{v}\circ \widetilde{\pi}_{v}$,\\
$\widetilde{\mathbf{M}}_{v}^{0}:=(\widetilde{a}_{v})^{-1}(v_{2},0)$,\\
$\mathbf{M}_{v}^{0}:=(a_{v})^{-1}(v_{2},0)$,\\
$\widetilde{\pi}_{v}^{0}:\widetilde{\mathbf{M}}_{v}^{0}\longrightarrow
\mathbf{M}_{v}^{0}$ the restriction of $\widetilde{\pi}_{v}$,\\
$\Sigma_{v}^{0}:=\mathbf{M}_{v}^{0}\cap\Sigma_{v}$,\\
$\Omega_{v}^{0}:=\mathbf{M}_{v}^{0}\cap\Omega_{v}$,\\
$\widetilde{\Sigma}_{v}^{0}:=\widetilde{\mathbf{M}}_{v}^{0}\cap
\widetilde{\Sigma}_{v}=\widetilde{\pi}_{v}^{-1}\Sigma_{v}^{0}$,\\
$\widetilde{\Omega}_{v}^{0}:=\widetilde{\mathbf{M}}_{v}^{0}\cap
\widetilde{\Omega}_{v}=\widetilde{\pi}_{v}^{-1}\Omega_{v}^{0}$.\\
When we will consider the case
\begin{equation}\label{vogrady}v=(2,0,-2)\end{equation}
(namely the Mukai vector used by O'Grady) we will generally use
the original notation $\widetilde{\mathcal{M}}$,
$\widetilde{\Sigma}$ and $\widetilde{\pi}$ to denote
$\widetilde{\mathbf{M}}_{v}^{0}$ and $\widetilde{\Sigma}^{0}_{v}$
and $\widetilde{\pi}_{v}^{0}$ respectively.

\end{notation}
\begin{rem}\label{notazioni fm} We will often implicitly use the identification
$\mathcal{J}\simeq\widehat{\mathcal{J}}$ given by means of
$\Theta$. In particular the Fourier-Mukai transform $\mathcal{FM}$
will be seen as the self-equivalence of the derived category of
coherent sheaves on $\mathcal{J}$ induced by the functor \[F
\mapsto q_{2*}(\mathcal{P}\otimes q_{1}^{*}F)\] where
$q_{1}:\mathcal{J}\times \mathcal{J}\rightarrow \mathcal{J}$ is
the projection on the ith-factor and
$\mathcal{P}:=q_{1}^{*}\mathcal{O}(\Theta)\otimes
q_{2}^{*}\mathcal{O}(\Theta)\otimes m^{*}\mathcal{O}(\Theta)^{\vee}$ is
the Poincar\'e line bundle
($m:\mathcal{J}\times\mathcal{J}\rightarrow \mathcal{J}$ is the
sum). Analogously the map induced in cohomology by the
Fourier-Mukai transform will be seen as an endomorphism of
$H^{2\bullet}$: precisely the endomorphism given by $\alpha\mapsto
q_{2*}(ch(\mathcal{P})\otimes q_{1}^{*}\alpha)$.\\ Finally, given
a sheaf $F$ on $\mathcal{J}$ we will denote by $c_{i}(F)$
($ch_{i}(F)$) the i-th Chern class (degree i component of the
Chern character) of F in both the cohomology ring and the Chow
ring: we will usually consider Chern classes (characters) as
cohomology classes and in the opposite case we will explicitly say
that we are referring to classes in the Chow ring.
\end{rem}
\section{Euler characteristic of $\widetilde{\mathcal{M}}$}

In this section we prove (Theorem \ref{caratteristica m tilde})
that the Euler characteristic of $\widetilde{\mathcal{M}}$ is
$1920$. By Corollary \ref{vettori espliciti}
$\widetilde{\mathbf{M}}_{(0,2\Theta,-2)}^{0}$ is a symplectic
desingularization of $\mathbf{M}_{(0,2\Theta,-2)}^{0}$. In
subsection 1 we study $\mathbf{M}_{(0,2\Theta,-2)}^{0}$ and
determine its Euler characteristic. Since we know the fibers of
the desigularization morphism (see Remark \ref{fibre ogradi}), we
can also deduce the Euler characteristic of
$\widetilde{\mathbf{M}}_{(0,2\Theta,-2)}^{0}$ (see
\ref{caratteristica m tilde02-2}) .\\ In subsection 2 we give an
explicit birational map between
$\widetilde{\mathbf{M}}_{(0,2\Theta,-2)}^{0}$ and
$\widetilde{\mathcal{M}}$. It turns out that
$\widetilde{\mathbf{M}}_{(0,2\Theta,-2)}^{0}$ is an irreducible
symplectic variety too. By a theorem due to Huybrechts the
existence of a birational map implies that
$\widetilde{\mathcal{M}}$ and
$\widetilde{\mathbf{M}}_{(0,2\Theta,-2)}^{0}$ are deformation
equivalent: hence their Euler characteristics are equal.
\subsection{Euler characteristic of $\widetilde{\mathbf{M}}_{(0,2\Theta,-2)}^{0}$ }

The analysis of $\mathbf{M}_{(0,2\Theta,-2)}^{0}$ is simplified by
the existence of the regular morphism
\begin{equation}\label{fitting} \Phi \colon \mathbf{M}_{(0,2\Theta,-2)}^{0}\rightarrow
|2\Theta|\end{equation} associating to an S-equivalence class of
sheaves the fitting subscheme of each of its representative.
\begin{rem}\label{fitting astra}
Recall that given a locally free presentation of a sheaf $F$,
\[F_{1}\stackrel{f}{\rightarrow} F_{0}\rightarrow F \rightarrow 0\]
the Fitting subscheme of $F$ is defined as the cokernel of the map
\[\wedge^{n}F_{1}\otimes \wedge^{n}F_{0}^{\vee}\rightarrow \mathcal{O}\]
($n$ being the rank of $F_{0}$) induced by $f$. In the case of
pure 1-dimensional sheaves on a smooth surface the construction of
the Fitting subscheme globalizes transforming flat families of
sheaves into flat families of 1-dimensional subschemes: thus it
induces regular maps between moduli spaces parametrizing pure
1-dimensional sheaves and Hilbert schemes parametrizing curves
(see \cite{LP93}). Moreover we can actually choose, as a locally
free presentation, a locally free resolution of $F$: this implies
that the fundamental cycle of the Fitting subscheme of $F$ is a
representative of $c_{1}(F)$ in the Chow ring of $\mathcal{J}$, in
particular, for $[F]\in
\widetilde{\mathbf{M}}_{(0,2\Theta,-2)}^{0}$, its Fitting
subscheme belongs to $|2\Theta|$.
\end{rem}
In order to study $\widetilde{\mathbf{M}}_{(0,2\Theta,-2)}^{0}$,
we now describe singularities occurring in curves in the linear
system $|2\Theta|$.

\begin{lem}\label{classificazione curvelem}
Let $C\in |2\Theta|$:\begin{enumerate}\item{ If
$C\cap\mathcal{J}[2]=\emptyset$  then either $C$ is smooth or
$C=\Theta_{x}+ \Theta_{-x}$ and $\Theta_{x}\cap \Theta_{-x}$
consists of 2 distinct points} \item{If $p\in\mathcal{J}[2]$
belongs to $C$ then  $C$ is singular at $p$: in this case either
$C$ is an irreducible  nodal curve smooth outside of
$\mathcal{J}[2]$ or $C=\Theta_{x}+ \Theta_{-x}$ with
$\Theta_{x}\cap \Theta_{-x}=p$ or $\Theta_{x}= \Theta_{-x}$}
\end{enumerate}
\end{lem}
\begin{proof}
Recall that the linear system $|2\Theta|$ induces a map $f\colon
\mathcal{J}\rightarrow \mathbb{P}^{3}$ whose image $Kum_{s}$ is
just the Kummer surface of $\mathcal{J}$, namely the quotient of
$\mathcal{J}$ by the involution $-1$: so it is a nodal surface
singular in $f(\mathcal{J}[2])$.\\
Given $p$ not belonging to $\mathcal{J}[2]$ there exists a unique
curve $C\in|2\Theta|$ singular in $p$: indeed since $f$ is \'etale
outside  $\mathcal{J}[2]$, $C$ is singular in $p$ if and only if
$f(C)$ is singular in $f(p)$ and there is a unique plane section
of $Kum_{s}$ singular in  $f(p)$, namely the one obtained
intersecting  $Kum_{s}$ with the plane tangent to  $Kum_{s}$ in
$p$. Since it is easily proved that there exists a curve of the
form $C=\Theta_{x}+ \Theta_{-x}$ (for all $x$ such a curve belong
to $|2\Theta|$)
singular in $p$, it is the unique one.\\
If $C\cap\mathcal{J}[2]=\emptyset$ and $C$ singular in $p$, then
it is singular also in $-p$: therefore $\Theta_{x}\cap
\Theta_{-x}$ consists of, at least, 2 points. Moreover, since
$\Theta^{2}=2$, if $\Theta_{x}\cap \Theta_{-x}$ contained a third
point we could conclude $\Theta_{x}= \Theta_{-x}$ thus
$x\in\mathcal{J}[2]$ and
$x\in C$. This proves item 1).\\
If $p\in\mathcal{J}[2]$ belongs to $C$ then $f(p)$ is a node of
$Kum_{s}$ and $f(C)=H\cap Kum_{s}$, $H$ being  a plane of
$\mathbb{P}^{3}$ passing through $f(p)$. The normal cone to
$Kum_{s}$ in $f(p)$ is the cone over a conic: if $H$ intersects
this cone in $2$ distinct lines, $f(C)$ has a node in $f(p)$ and
$C$, being a $[2:1]$ covering ramified in $f(p)$, has a node in
$p$. The planes containing $f(p)$ and  not intersecting the normal
cone in $2$ distinct lines are parametrized by a conic: thus the
locus of $|2\Theta|$ parametrizing curves passing through $p$ and
not having a node in $p$ is contained in a conic. On the other
hand the curves of the form $\Theta_{p+x}+\Theta_{p-x}$ with $x\in
\Theta$ never have a node in $p$: in fact,  for $x\in \Theta$, an
easy computation on $Pic(C^{0})$ shows that either
$\Theta_{p+x}\cap\Theta_{p-x}=p$ or $x=-x$ and
$\Theta_{p+x}=\Theta_{p-x}$. Since, as already showed, any $C\in
|2\Theta|$ singular outside $\mathcal{J}[2]$ is of the form
$\Theta_{x}+\Theta_{-x}$, item (2) follows.

\end{proof}
Since $NS(\mathcal{J})=\mathbb{Z}\Theta$ the geometric genus of
curve of $\mathcal{J}$ is at least $2$. Since a curve $C$ in
$|2\Theta|$ has arithmetic genus $5$, it can contain at most $3$
nodes. Using this remark the previous Lemma can be reformulated as
follows.
\begin{prop}\label{classificazione curve} If
$NS(\mathcal{J})=\mathbb{Z}\Theta$, the stratification of
$|2\Theta|$ by the analytic type of singularity is the following:
\begin{itemize}
\item{Stratum S: the locus parametrizing smooth curves of genus
5.} \item{Stratum N(1): the locus parametrizing nodal irreducible
curves  singular  in a unique 2-torsion point.} \item{Stratum
N(2): the locus parametrizing nodal irreducible curves singular
only in 2 distinct  2-torsion points.} \item{Stratum N(3): the
locus parametrizing nodal irreducible curves singular only in 3
distinct points of 2-torsion.}\item{Stratum R(1): the locus
parametrizing reducible curves with nodal singularity (they are
the curves of the form $\Theta_{x}\cup\Theta_{-x}$ with 2 singular
points, namely $\Theta_{x}\cap\Theta_{-x}$ consists of 2 distinct
points).} \item{Stratum R(2): the locus parametrizing reducible
curves with a unique singular point (they are the curves of the
form $\Theta_{x}\cup\Theta_{-x}$ with a unique (non nodal)
singular point, namely $\Theta_{x}\cap\Theta_{-x}$ consists of a
unique point belonging to $\mathcal{J}[2]$).} \item{Stratum D: the
locus parametrizing non reduced curves (they are the curves of the
form $2\Theta_{x}$ with x a point of 2-torsion).}
\end{itemize}
\end{prop}
In order to determine the Euler characteristic of
$\widetilde{\mathcal{M}}$ we need to compute the Euler
characteristic of the fibers of $\Phi$ and establish their
dimension.
\begin{prop}\label{caratteristica fibre}
For any $C\in |2\Theta|$ the dimension of $\Phi^{-1}(C)$ is 3.\\
If $C\in S\cup N(1)\cup N(2)\cup R(1) \cup R(2)$ then
$\chi(\Phi^{-1}(C))=0$.\\
If $C\in N(3)$ then $\chi(\Phi^{-1}(C))=4$.\\
If $C\in D$ then $\chi(\Phi^{-1}(C))=20$.
\end{prop}
\begin{proof}
In this proof we will denote by $M_{C}$, the locus of
$\mathbf{M}_{(0,2\Theta,-2)}$ parametrizing sheaves whose Fitting
subscheme is $C$: in particular $\Phi^{-1}(C)=M_{C}\cap
\mathbf{M}_{(0,2\Theta,-2)}^{0}$.\\
If $C\in S\cup N(1)\cup N(2)\cup N(3)$, then points of $M_{C}$
correspond to isomorphism classes of rank-1, torsion-free sheaves
on $C$ whose Euler characteristic is $-2$. It is known (see
\cite{Be99}) that such a sheaf is either a degree-2 line bundle
on $C$ or the push-forward of a degree-(2-r) line bundle from a
partial
normalization desingularizing exactly r nodes of $C$.\\
Since line bundles having fixed degree on a nodal curve with n
nodes are parametrized by a $(\mathbb{C}^{*})^{r}-bundle$ on the
Jacobian of its normalization (see \cite{HM98}), $M_{C}$ can be
stratified in $(\mathbb{C}^{*})^{r}-bundle$ over the Jacobian
$J(\widetilde{C})$ of
the normalization $\widetilde{C}$.\\
Letting $C'$ be a partial normalization of $C$ we now determine
the intersection of the stratum $U(C')$ parametrizing push-forward
of line bundle on $C'$ with $\Phi^{-1}(C)$. The restriction
$r:U(C')\rightarrow \mathcal{J}$ of the map $a_{(0,2\Theta,-2)}$
to $U(C')$ associates to a sheaf $F$, the point $\sum
n_{i}p_{i}+\sum q_{k}$ where $q_{k}\in\mathcal{J}[2]\cap C$ are
points having 2 distinct inverse images on $C'$ and $\sum
n_{i}p_{i}$ is the push-forward on $\mathcal{J}$ of a
representative, in the Chow
ring of  $\widetilde{C}$, of $c_{1}$ of the pull-back of $F$.\\
It follows that $r$ descends to a map from $J(\widetilde{C})$:
moreover this last map is easily checked to be identified with the
map
\begin{equation}
a\colon J(\widetilde{C})\rightarrow \mathcal{J}
\end{equation}induced on the Albanese varieties by the morphism
$\widetilde{C}\rightarrow J$ obtained composing the normalization
morphism of $C$ with the inclusion of $C$ in $\mathcal{J}$.
$U(C')$ is then a $(\mathbb{C}^{*})^{r}-bundle$ over the fiber of
$a$: if $C \in S\cup N(1)\cup N(2)$ then
$dim(J(\widetilde{C}))>dim(\mathcal{J})$ and the Euler
characteristic of a fiber of $a$ is $0$: it follows that
$\chi(U(C'))=0$ too.\\ Since $\Phi^{-1}(C)$ is stratified by
subvarieties of the form $U(C')$, by the additivity of the Euler
characteristic, we get $\chi(\Phi^{-1}(C))=0$.\\
If $C\in N(3)$ the strata of $\Phi^{-1}(C)$ parametrizing sheaves
which are not push-forward of line bundle from $\widetilde{C}$,
are bundles with fiber $(\mathbb{C}^{*})^{r}$ with $r>0$: their
Euler characteristics are still 0. Using again the additivity of
the Euler characteristic we get that $\chi(\Phi^{-1}(C))$ equals
the Euler characteristic of the  stratum  parametrizing sheaves
which are push-forward of line bundle from $\widetilde{C}$.\\
This stratum is, as already explained, in bijective correspondence
with $a^{-1}(0)$. Since $C\in N(3)$, $J(\widetilde{C})$ is an
abelian surface and since $NS(\mathcal{J})=\mathbb{Z}\Theta$,
$\mathcal{J}$ doesn't contain elliptic curves: therefore $a\colon
J(\widetilde{C})\rightarrow \mathcal{J}$ is an \'etale covering
and $J(\widetilde{C})$ has the same Hodge structure over
$\mathbb{Q}$ of $\mathcal{J}$, in particular letting
$\Theta_{\widetilde{C}}$ be the theta divisor on
$J(\widetilde{C})$ we get
$NS(J(\widetilde{C}))=\mathbb{Z}\Theta_{\widetilde{C}}$. The
degree of $a$ is easily computed to be 4: in fact
$a^{*}(\Theta)\cap \Theta_{\widetilde{C}}=C\cap
\Theta=4=2\Theta_{\widetilde{C}}^{2}$, hence
$a^{*}(\Theta)=2\Theta_{\widetilde{C}}$ and for $C\in N(3)$
\begin{equation}\chi(\Phi^{-1}(C))=deg(a)=
\frac{(2\Theta_{\widetilde{C}})^{2}}{(\Theta)^{2}}=4.
\end{equation}
Moreover, using again that $\mathcal{J}$ doesn't contain elliptic
curves we get that $a$ is always surjective, hence the fiber of
$a$ has dimension equal to the geometric genus of $C$ minus 2: it
follows from the given description of the stratification  of
$\Phi^{-1}(C)$, that for any $C\in S\cup N(1)\cup N(2)\cup N(3)$
the stratum of $\Phi^{-1}(C)$ parametrizing line bundles on $C$
has dimension 3 and its complement has lower dimension.\\
We now consider the case $C\in  R(1) \cup R(2)$. In this case
$C=\Theta_{x}+\Theta_{-x}$. We will denote by
$i_{x}:C^{0}\rightarrow \mathcal{J}$ and $i_{-x}:C^{0}\rightarrow
\mathcal{J}$ the embeddings of $C^{0}$ on $\Theta_{x}$ and
$\Theta_{-x}$ respectively.\\
Since the Fitting subscheme of a sheaf contains the subscheme
defined by the annihilator of the sheaf (see \cite{Ei95}), for
any sheaf $F$ such that $[F]\in \Phi^{-1}(C)$ is the push-forward
via the inclusion $i:C\rightarrow \mathcal{J}$ of a sheaf $F_{C}$
on $C$: moreover since $c_{1}(F)=2\Theta$, the restriction of
$F_{C}$ to each of the component of $C$ is a rank-1 sheaf.\\
By the description of strictly semistable sheaves given in the
proof of \ref{desingk3}, points the locus of $M_{C}$ parametrizing
strictly semistable sheaves are in bijective correspondence with
isomorphism classes of sheaves of the form $G_{1}\oplus G_{2}$,
where $G_{i}$ are stable sheaves whose Mukai vector is
$(0,c_{1}(\Theta),-1)$, namely $G_{1}=(i_{x})_{*}L_{1}$ and
$G_{2}=(i_{-x})_{*}L_{2}$ where $L_{i}$ are degree-0 line bundles
on $C^{0}$. With these notations, fixed $L_{1}\in Pic^{0}(C^{0})$,
since the map associating to each line bundle $L$ on $C^{0}$ the
point $\sum n_{i}p_{i}\in\mathcal{J}$ (where $\sum n_{i}p_{i}$ is
a representative of $c_{2}((i_{x})_{*}L)$ in the Chow ring) is
obviously an isomorphism, we get that there exists a unique
$L_{2}$ such that $[(i_{x})_{*}L_{1}\oplus(i_{-x})_{*}L_{2}]\in
\Phi^{-1}(C)$. It follows that the strictly semistable locus
$\Phi^{-1}(C)^{ss}\subset\Phi^{-1}(C)$ is isomorphic to
$\mathcal{J}$.\\
To study the stable locus $\Phi^{-1}(C)^{s}$, for any $F$ having
$C$ as Fitting subscheme, we denote by $L_{1}$ and $L_{2}$  the
torsion free parts of $i_{x}^{*}(F)$ and $i_{-x}^{*}(F)$. There
are natural surjective maps from $F$ to $(i_{x})_{*}L_{1}$ and to
$(i_{-x})_{*}L_{2}$: their direct sum is the first map in the
following exact sequence
\begin{equation}\label{quisotto}
0\rightarrow F \stackrel{\alpha}{\rightarrow}
(i_{x})_{*}L_{1}\oplus (i_{-x})_{*}L_{2}
\stackrel{\beta}{\rightarrow} Q \rightarrow 0.
\end{equation}
Since both the components of $\alpha$ are surjective $Q$ is a
quotient of both $(i_{x})_{*}L_{1}$ and $(i_{x})_{*}L_{2}$  and
therefore a quotient of $\mathcal{O}_{\Theta_{x}\cap\Theta_{-x}}$.
If $F$ is stable then $deg(L_{i})\ge 1$, and using
(\ref{quisotto}) to compute Chern classes of $F$ we find
$deg(L_{1})+deg(L_{2})=lenght(Q)$: is then easy to check that
$lenght(Q)=2$ $(Q=\mathcal{O}_{\Theta_{x}\cap\Theta_{-x}})$ and
$deg(L_{2})=deg(L_{1})=1$. It follows that if $F$ is  stable it is
the kernel of a map $\beta\colon (i_{x})_{*}L_{1}\oplus
(i_{x})_{*}L_{2}\rightarrow
\mathcal{O}_{\Theta_{x}\cap\Theta_{-x}}$, where $L_{i}$ is a
degree-1 line bundle and the restriction of $\beta$ to each
summand is already surjective.\\
On the other hand any such a kernel is easily seen to be a stable
sheaf, and, given 2  kernels obtained in this way, they are
isomorphic if and only if they differ by an automorphism of
$(i_{-x})_{*}L_{1}\oplus (i_{x})_{*}L_{2}$.\\
Fixed $L_{1}$ and $L_{2}$ it can be seen the following: if $C\in
R(1)$ the kernels obtained are parametrized by
$\mathbb{C}^{*}\times\mathbb{C}^{*}$ and  the isomorphism classes
of sheaves simply by $\mathbb{C}^{*}$, if $C\in R(2)$ the kernels
obtained are parametrized by $\mathbb{C}^{*}\times\mathbb{C}$ and
the isomorphism classes of sheaves simply by $\mathbb{C}$.\\
As in the strictly semistable case we can see that, fixed $L_{1}$
there exists a unique $L_{2}$ giving kernels belonging to
$\Phi^{-1}(C)$.\\
It follows that, for $C\in R(1)\cup R(2)$,  $\Phi^{-1}(C)^{s}$ is
either a $\mathbb{C}^{*}$-bundle or a $\mathbb{C}$-bundle over
$\mathcal{J}$. Since $\Phi^{-1}(C)$ is the disjoint union of
$\Phi^{-1}(C)^{s}$ and $\mathcal{J}$, we get
$\chi(\Phi^{-1}(C))=0$ and $dim(\Phi^{-1}(C))=3$.\\
It remains to consider the case $C\in D$. To simplify the notation
we deal explicitly only with the case $C=2\Theta$, but the same proof
works with $C$ replaced by any double curve in $|2\Theta|$.\\ In
this case the subscheme defined by the annihilator of a sheaf $F$
such that $[F]\in\Phi^{-1}(C)$, being a pure 1-dimensional scheme
contained in $C$ is either $\Theta$ or the double curve $C$.\\
If the annihilator of $F$ (such that $[F]\in\Phi^{-1}(C)$) is the
ideal $I_{\Theta}$ of $\Theta$ then $F=(i_{0})_{*}V$, where $V$ is
a semistable rank-2  vector bundle with trivial determinant and
$i_{0}\colon C^{0}\rightarrow \mathcal{J}$ is the imbedding of
$C^{0}$ on $\Theta$. The locus of $\Phi^{-1}(C)$ parametrizing
sheaves of this form is then isomorphic to the moduli space of
rank-2  semistable  vector bundles with fixed determinant on a
genus 2 curve: Narashiman and Ramanan proved in \cite{NR69} that
it isomorphic
to $\mathbb{P}^{3}$.\\
Since any polystable sheaf $F$ such that $[F]\in\Phi^{-1}(C)$
is the push-forward of a rank-2 vector bundle from $C$,
points of the complement $U$ of this $\mathbb{P}^{3}$ are in [1:1]
correspondence with isomorphism classes of stable sheaves
annihilated by the ideal $I_{C}$ of $C$.\\
Letting $F$ be such a sheaf we want to prove that it fits in an exact sequence
\begin{equation}\label{carattezzazione fasci doppi}
0\rightarrow (i_{0})_{*}(K^{\vee}\otimes L) \rightarrow  F
\rightarrow  (i_{0})_{*} L\rightarrow 0.
\end{equation}
where $K$ is the canonical bundle and $L$ is a degree-1 line
bundle such that
$L^{\otimes2}\otimes K^{\vee}=\mathcal{O}_{C^{0}}$.\\
Making the Tensor product of $F$ with the exact sequence of
sheaves on $\mathcal{J}$ defining the double structure of $C$
\begin{equation}
0\rightarrow (i_{0})_{*}K^{\vee}\rightarrow \mathcal{O}_{C}
\rightarrow  \mathcal{O}_{2\Theta}\rightarrow 0.
\end{equation}
we get the exact sequence
\begin{equation}
(i_{0})_{*}(K^{\vee}\otimes L \oplus T)\rightarrow  F
\rightarrow  (i_{0})_{*}(L\oplus T)\rightarrow 0.
\end{equation}
where  $L$ is a vector bundle and $T$ a torsion sheaf on $C^{0}$.
Since $c_{1}(F)=2\Theta$  and $F$ is not the push-forward of a
sheaf from $C^{0}$, the rank of $L$ is $1$. Since $F$ is stable
and $F$ surjects on $(i_{0})_{*} L$, we have $deg(L)>0$. On the
other hand since $F$ is pure 1-dimensional the kernel of the first
map in the last exact sequence is just $T$. Since $ch_{2}(F)=-2$
we get $2deg(L)+length(T)=2$: this implies $deg(L)=1$ and $T=0$,
so proving the existence of the exact sequence
(\ref{carattezzazione fasci doppi}) ($L^{\otimes2}\otimes
K^{\vee}=\mathcal{O}_{C^{0}}$ is required to have $[F]\in
\Phi^{-1}(C)$). Moreover it is easily seen that for any sheaf $F$,
having the ideal $I_{C}$ as its annihilator and fitting in such an
exact sequence a pure 1-dimensional sheaf and the subsheaf
$(i_{0})_{*}(K^{\vee}\otimes L)$ is just the subsheaf annihilated
by the ideal $I_{\Theta}$ of $\Theta$: since stabilty can be
checked using only injections of push-forwards of line bundle on
$C^{0}$ and these are always annihilated by $I_{\Theta}$ it
follows that any $F$ under our condition is stable.\\
We can conclude that for any $L$ satisfying $L^{\otimes2}\otimes
K^{\vee}=\mathcal{O}_{C^{0}}$ there is a locus $U_{L}\subset U$
whose points are in bijective correspondence with isomorphism
classes of extensions of the form (\ref{carattezzazione fasci
doppi}) with $F$ not annihilated by $I_{\Theta}$. Using
Hirrzebruch-Riemann-Roch and Serre duality for extensions (see
\cite{HL97}) we get
$dim(Ext^{1}((i_{0})_{*}L,(i_{0})_{*}(K^{\vee}\otimes L))=4$ and
the extensions being push-forward of vector bundles can be
identified using the spectral sequence associated to the
composition
 $H^{0}\circ
\mathcal{H}om(i_{*}L,\cdot)=Hom((i_{0})_{*}L,)$: it produces the
following short exact sequence
\[0\rightarrow
H^{1}(\mathcal{H}om((i_{0})_{*}L,(i_{0})_{*}(K^{\vee}\otimes L)))
\rightarrow\]\[ Ext^{1}((i_{0})_{*}L,(i_{0})_{*}(K^{\vee}\otimes
L))\rightarrow
H^{0}(\mathcal{E}xt^{1}((i_{0})_{*}L,(i_{0})_{*}(K^{\vee}\otimes
L))) \rightarrow 0 .\] where the first term is isomorphic to
$\mathbb{C}^{3}$ and parametrizes just the extensions coming from
extensions of line bundles on $C^{0}$.\\ Since for any extension
belonging to $Ext^{1}((i_{0})_{*}L,(i_{0})_{*}(K^{\vee}\otimes L))
\setminus
H^{1}(\mathcal{H}om((i_{0})_{*}L,(i_{0})_{*}(K^{\vee}\otimes L)))$
we have $End((i_{0})_{*}(K^{\vee}\otimes
L))=End((i_{0})_{*}L)=End(F)=\mathbb{C}$, 2 such extensions have
middle terms isomorphic if and only if they differ by a scalar
multiplication: it follows that $U_{L}$ is in bijective
correspondence with
$\mathbb{C}^{3}$.\\
Since there are 16 line line bundles on $C^{0}$ satisfying
$L^{\otimes2}\otimes K^{\vee}=\mathcal{O}_{C^{0}}$ and the
respective $U_{L}'s$ are easily seen to be disjoint we get that
$\Phi^{-1}(C)$ is the disjoint union of $\mathbb{P}^{3}$ with 16
3-dimensional affine spaces: therefore $dim(\Phi^{-1}(C))=3$ and
$\chi(\Phi^{-1}(C))=20$.
\end{proof}
Before calculating the Euler characteristic of
$\mathbf{M}_{(0,2\Theta,-2)}^{0}$ we establish its irreducibility
as a corollary of the previous proposition.
\begin{cor} \label{irriducibilita m02-2}
$\mathbf{M}_{(0,2\Theta,-2)}^{0}$ and
$\mathbf{M}_{(0,2\Theta,-2)}$ are reduced irreducible.
\end{cor}
\begin{proof}Since $a_{(0,2\Theta,-2)}$ makes $\mathbf{M}_{(0,2\Theta,-2)}$
a fibration over $\mathcal{\widehat{J}}\times \mathcal{J}$ with
fiber $\mathbf{M}_{(0,2\Theta,-2)}^{0}$, the second statement
follows
from the first.\\
Let's prove  that $\mathbf{M}_{(0,2\Theta,-2)}$ is reduced
irreducible.\\
 Since by Proposition \ref{desingk3}
$\mathbf{M}_{(0,2\Theta,-2)}^{0}$ is reduced purely 6-dimensional,
it is enough to prove that there exists an irreducible open
subvariety $U\subset\mathbf{M}_{(0,2\Theta,-2)}^{0}$ whose
complement has dimension at most $5$.\\ We set $U:=\Phi^{-1}(S)$.
By the previous proposition we have
$dim(\mathbf{M}_{(0,2\Theta,-2)}^{0}\setminus U)\ge 5$. To show
the irreducibility of $U$ recall that for $C\in S$ the fiber
$\Phi^{-1}(C)$ is identified to the kernel of the map
\[a:
J(C):=\frac{H^{1}(\Omega_{C})^{\vee}}{H_{1}(C,\mathbb{Z})}\rightarrow
\frac{H^{1}(\Omega_{\mathcal{J}})^{\vee}}{H_{1}(\mathcal{J},\mathbb{Z})}\]
induced on the Albanese varieties by the embedding of $C$ in
$\mathcal{J}$. The Kernel of $a$ is an irreducible torus because
$C$ is an ample divisor and hence, by the Hyperplane section
theorem, $H_{1}(C,\mathbb{Z})$ surjects on
$H_{1}(\mathcal{J},\mathbb{Z})$. $U$ is finally irreducible being
a bundle with irreducible fiber and base.
\end{proof}
We can now compute the Euler characteristic of
$\mathbf{M}_{(0,2\Theta,-2)}^{0}$.
\begin{prop}\label{caratteristica moduli}
$\chi(\mathbf{M}_{(0,2\Theta,-2)}^{0})=1280$.
\end{prop}
\begin{proof}
By Proposition \ref{classificazione curve} and the additivity of
Euler characteristic we get
$\chi(\mathbf{M}_{(0,2\Theta,-2)}^{0})=\chi(\Phi^{-1}(S\cup N(1)
\cup N(2) \cup R(1) \cup R(2))+\chi(D\cup N(3))$ and since the
Euler characteristic of any fiber in $\Phi^{-1}(S\cup N(1) \cup
N(2) \cup R(1) \cup R(2)$ is $0$, we obtain $\chi(\Phi^{-1}(S\cup
N(1) \cup N(2) \cup R(1) \cup R(2))=0$ (see \cite{Be99}). Letting
$|N(3)|$ be the cardinality of $N(3)$ and using again Proposition
\ref{caratteristica fibre} we conclude
\begin{equation}\label{111}\chi(\mathbf{M}_{(0,2\Theta,-2)}^{0})
=\chi(\Phi^{-1}(D))+\chi(\Phi^{-1}(n(3)))=16\cdot20+4|N(3)|.\end{equation}
$|N(3)|$ is computed as follows. A triple of singular points of
$Kum_{s}$ defines a curve in $N(3)$ if and only if it is not
included in a double curve: thus
\begin{equation}\label{222} |N(3)|=\binom{16}{3}-16\binom{6}{3}
=240,\end{equation}
$\binom{6}{3}$ are the triples of singular points included in a
double curve.\\
Formulas (\ref{111}) and  (\ref{222}) and Proposition
\ref{caratteristica fibre} imply the result.
\end{proof}
Since we know  the birational modification needed to obtain
$\widetilde{\mathbf{M}}_{(0,2\Theta,-2)}^{0}$ from
$\mathbf{M}_{(0,2\Theta,-2)}^{0}$ it is now possible to compute
$\chi(\widetilde{\mathbf{M}}_{(0,2\Theta,-2)}^{0})$.
\begin{thm}\label{caratteristica m tilde02-2}
$\chi(\widetilde{\mathbf{M}}_{(0,2\Theta,-2)}^{0})=1920$.
\end{thm}
\begin{proof}
Since $\widetilde{\mathbf{M}}_{(0,2\Theta,-2)}^{0}$ is the
disjoint union of the stable locus  of
$\mathbf{M}_{(0,2\Theta,-2)}^{0}$ and
$\widetilde{\Sigma}_{(0,2\Theta,-2)}^{0}$ and by the additivity of
$\chi$ we have
\[\chi(\widetilde{\mathbf{M}}_{(0,2\Theta,-2)}^{0})=
\chi(\mathbf{M}_{(0,2\Theta,-2)}^{0})-
\chi(\Sigma_{(0,2\Theta,-2)}^{0})+
\chi(\widetilde{\Sigma}_{(0,2\Theta,-2)}^{0}).\]  Since the map
$\mathcal{J}\times \mathcal{\widehat{J}}\rightarrow
\Sigma_{(0,2\Theta,-2)}^{0}$, associating to $(x,L)$ the the
s-equivalence class of $i_{x*}L\oplus i_{-x*}L$, is surjective and
[2:1] outside  the 256 2-torsion points of $\mathcal{J}\times
\mathcal{\widehat{J}}$ we find
\[\chi(\Sigma_{(0,2\Theta,-2)}^{0})=
\frac{\chi(\mathcal{J}\times \mathcal{\widehat{J}})+256}{2}=128.\]
To compute $\chi(\widetilde{\Sigma}_{(0,2\Theta,-2)}^{0}) $ recall
( see remark \ref{fibre ogradi}) that the restriction to
$\widetilde{\Sigma}_{(0,2\Theta,-2)}^{0} $ of the
desingularization map is a $\mathbb{P}^{1}$ bundle outside the 256
points of $\Omega_{(0,2\Theta,-2)}^{0}$ where the fibers are
smooth 3-dimensional quadrics. The Euler characteristic of the
$\mathbb{P}^{1}$ bundle is $-256$ and the Euler characteristic of
the 3-dimensional quadric is 4: the final result is
\[\chi(\widetilde{\mathbf{M}}_{(0,2\Theta,-2)}^{0})=
1280-128-256+256\cdot4=1920.\]
\end{proof}
\subsection{$\widetilde{\mathbf{M}}_{(0,2\Theta,-2)}^{0}$
is a birational model of $\widetilde{\mathcal{M}}$}

In this section we show that
$\chi(\widetilde{\mathbf{M}}_{(0,2\Theta,-2)}^{0})=
\chi(\widetilde{\mathcal{M}})$.\\
This equation is a consequence of the following proposition.
\begin{prop}\label{v02-2}
There exists a birational map
$b:\widetilde{\mathcal{M}}\dashrightarrow
\widetilde{\mathbf{M}}_{(0,2\Theta,2)}^{0}$.\\
Moreover, letting $b^{*}$ be the pull-back of divisors,
$b^{*}\widetilde{\Sigma}_{(0,2\Theta,-2)}^{0}=\widetilde{\Sigma}$.
\end{prop}
\begin{proof}
Since $\mathbf{NS}(\mathcal{J})=\mathbb{Z}\Theta$ the
tensorization by $\mathcal{O}(\Theta)$  doesn't change stability
and semistability, hence it  induces an isomorphism
\[t:\mathbf{M}_{(2,0,-2)}\rightarrow\mathbf{M}_{(2,2\Theta,0)}.\]
An easy calculation shows that this isomorphism  sends
$\Sigma_{(2,0,-2)}$ to $\Sigma_{(2,2\Theta,0)}$ and
$\mathbf{M}_{(2,0,-2)}^{0}$ to $\mathbf{M}_{(2,2\Theta,0)}^{0}$.\\
Recalling that, for $v$  as in Proposition \ref{desing abel} the
blow up
$Bl_{\Sigma^{0}_{v}\setminus\Omega^{0}_{v}}\mathbf{M}^{0}_{v}\setminus\Omega^{0}_{v}$
is an open subset of $\widetilde{\mathbf{M}}^{0}_{v}$ and its
complement has codimension bigger than $1$ (see remark
\ref{costruzione og rem}), we get that $t$ induces a birational
map (actually biregular)
$\widetilde{t}:\widetilde{\mathcal{M}}\rightarrow
\widetilde{\mathbf{M}}_{(2,2\Theta,0)}^{0}$ such that
$\widetilde{t}^{*}\widetilde{\Sigma}_{(2,2\Theta,0)}^{0}=\widetilde{\Sigma}$.\\
It remains to prove that there exists a birational map
$\widetilde{fm}:
\widetilde{\mathbf{M}}_{(2,2\Theta,0)}^{0}\dashrightarrow
\widetilde{\mathbf{M}}_{(0,2\Theta,-2)}^{0}$ such that
\[\widetilde{fm}^{*}\widetilde{\Sigma}_{(0,2\Theta,-2)}^{0}=
\widetilde{\Sigma}_{(2,2\Theta,0)}^{0}.\] We have denoted this map
by $\widetilde{fm}$ because it is induced by the
Fourier-Mukai transform (see remark \ref{notazioni fm})
\[\mathcal{FM}:D\mathcal{C}oh(\mathcal{J})\rightarrow
D\mathcal{C}oh(\mathcal{J}).\]\\
We need the following lemma.
\begin{lem} Let $F$ be a strictly semistable sheaf, having Mukai's vector
$(2,2\Theta,0)$ (or $(0,2\Theta,-2)$), then $F$ verifies the
W.I.T. (with index 1), and moreover $\mathcal{FM}(F)(-1)$ is a
strictly semistable sheaf having Mukai's vector $(0,2\Theta,-2)$
(or $(2,2\Theta,0)$).
\end{lem}
\begin{proof}
We only deal with the first case, the second being completely
analogous. By the description of strictly semistable sheaves
(sequence (3.3) in the proof of Proposition \ref{desingk3} ) it is
enough to verify  that a sheaf of the form $I_{x}\otimes
\mathcal{O}(\Theta_{y})$ verifies the W.I.T. (with index 1)  and
its Fourier-Mukai transform is a sheaf of the form $(i_{z})_{*}L$
where $i_{z}:C^{0}\rightarrow\mathcal{J}$ is the embedding on
$\Theta_{z}$ and $L$ is a degree-0 line bundle.\\
Indeed, since both $\mathbb{C}_{x}$ and $\mathcal{O}(\Theta_{y})$
satisfy the weak index theorem with index $0$ and their
Fourier-Mukai transform are $\mathcal{O}(\Theta_{x}-\Theta)$ and
$\mathcal{O}(\Theta_{y})^{\vee}$ respectively (see Theorem 3.13
\cite{Mu81}), the short exact sequence
\[0\rightarrow I_{x}\otimes\mathcal{O}(\Theta_{y})\rightarrow \mathcal{O}(\Theta_{y})
\rightarrow \mathbb{C}_{x}\rightarrow 0\;(*)\] induces the long
exact sequence
\[0\rightarrow q_{2*}(\mathcal{P}\otimes q_{1}^{*}(I_{x}\otimes\mathcal{O}(\Theta_{y})))
\rightarrow \mathcal{O}(\Theta_{y})^{\vee}\rightarrow
\mathcal{O}(\Theta_{x}-\Theta)\rightarrow
R^{1}q_{2*}(\mathcal{P}\otimes
q_{1}^{*}(I_{x}\otimes\mathcal{O}(\Theta_{y})))\rightarrow 0\] and
since the middle map of this sequence cannot be zero, the first
term is $0$ and the last one is just
$\mathcal{FM}(I_{x}\otimes\mathcal{O}(\Theta_{y}))(-1)$
and has the stated form.\\
\end{proof}
By general results on Fourier-Mukai transform (see \cite{Mu81}
and \cite{Mu87}), and since $\mathbf{M}_{(2,2\Theta,0)}$ and
$\mathbf{M}_{(0,2\Theta,-2)}$ are both reduced irreducible (see
\ref{irriducibilita m02-2}) this lemma implies that there exists a
birational map $fm: \mathbf{M}_{(2,2\Theta,0)}\rightarrow
\mathbf{M}_{(0,2\Theta,-2)}$ restricting to an isomorphism on a
neighborhood of $\Sigma_{(0,2\Theta,-2)}$ and such that
$fm(\Sigma_{(2,2\Theta,0)})=\Sigma_{(0,2\Theta,-2)})$.\\
Recalling how O'Grady's desingularizations are obtained, we get a
lift \[\overline{fm}:
\widetilde{\mathbf{M}}_{(2,2\Theta,0)}\rightarrow
\widetilde{\mathbf{M}}_{(0,2\Theta,-2)}\] such that
$\overline{fm}^{*}\Sigma_{(0,2\Theta,-2)}=\Sigma_{(2,2\Theta,0)}$.\\
To complete the proof of this proposition it remains to show that
$\overline{fm}$ sends  fibers of $a_{(2,2\Theta,0)}:
\widetilde{\mathbf{M}}_{(2,2\Theta,0)}\rightarrow\mathcal{J}\times\widehat{\mathcal{J}}$
birationally to fibers of $a_{(0,2\Theta,-2)}:
\widetilde{\mathbf{M}}_{(0,2\Theta,-2)}\rightarrow\mathcal{J}\times\widehat{\mathcal{J}}$.\\
Since the fibers of $a_{(2,2\Theta,0)}$ are isomorphic to
$\widetilde{\mathcal{M}}$, their fundamental group is trivial:
hence $a_{(2,2\Theta,0)}$ is identified  with the Albanese map of
$\widetilde{\mathbf{M}}_{(2,2\Theta,0)}$. On the other hand the
rational map $a_{(2,2\Theta,0)}\circ\overline{fm}:
\widetilde{\mathbf{M}}_{(2,2\Theta,0)}\rightarrow\mathcal{J}\times\widehat{\mathcal{J}}$,
being a map to an abelian variety, extends. By the universal
property of the Albanese map there exists a morphism $g$ making
commutative the following diagram:\[\CD
\widetilde{\mathbf{M}}_{(2,2\Theta,0)} @>\overline{fm}>>
\widetilde{\mathbf{M}}_{(0,2\Theta,-2)}\\
@Va_{(2,2\Theta,0)}VV  @Va_{(0,2\Theta,-2)}VV \\
\mathcal{J}\times\widehat{\mathcal{J}}@>g>> \mathcal{J}\times\widehat{\mathcal{J}}\\
\endCD
\] Thus $\overline{fm}$ sends fibers of $a_{(2,2\Theta,0)}$ to
fibers of $a_{(0,2\Theta,-2)}$. On the other hand we already
proved (Corollary \ref{irriducibilita m02-2}) that fibers of
$a_{(0,2\Theta,-2)}$ are irreducible. It follows that
$\overline{fm}$ induces, by restriction to the central fibers, a
map
$\widetilde{fm}:\widetilde{\mathbf{M}}^{0}_{(2,2\Theta,0)}\rightarrow
\widetilde{\mathbf{M}}^{0}_{(0,2\Theta,-2)}$ such that
$\widetilde{fm}^{*}(\Sigma^{0}_{(0,2\Theta,-2)})=\Sigma^{0}_{(2,2\Theta,0)}$.
The map $b:=\widetilde{fm}\circ \widetilde{t}$ verifies the thesis
of the  proposition.
\end{proof}
\begin{thm}\label{caratteristica m tilde}
$\chi(\widetilde{\mathbf{M}}^{0}_{(2,0,-2)})=1920$
\end{thm}
\begin{proof}
By Proposition \ref{desing abel} $\widetilde{\mathcal{M}}$ is
symplectic and projective. By the previous proposition it is
birational to an irreducible symplectic variety: it follows that
$\widetilde{\mathbf{M}}^{0}_{(0,2\Theta,-2)}$ is an irreducible
symplectic variety too.\\
By a theorem due to Huybrechts (\cite{Hu97}) birational
irreducible symplectic varieties are deformation equivalent:
therefore $\chi(\widetilde{\mathcal{M}})=
\chi(\widetilde{\mathbf{M}}^{0}_{(0,2\Theta,-2)})$ and the last is
$1920$ as shown in Theorem \ref{caratteristica m tilde02-2}.
\end{proof}
\section{The Beauville form of $\widetilde{\mathcal{M}}$}

In this section  we determine the Beauville form and the Fujiki
constant of $\widetilde{\mathcal{M}}$ (see Theorem \ref{mainthm}).
Before going on we recall the Theorem due to Beauville and Fujiki (see
\cite{Be83,Fu87}) that defines the the Beauville form and the Fujiki
constant of an irreducible symplectic variety.

\begin{thm}\label{teorema def forma beauville}
Let $X$ be a 2n dimensional irreducible symplectic manifold. There
exist a unique indivisible bilinear integral symmetric form
$B_{X}\in S^{2}(H^{2}(X,\mathbb{Z}))^{*}$ and a unique positive
constant $c_{X}\in\mathbb{Q}$ such that for any $\alpha\in
H^{2}(X,\mathbb{Z})$
\begin{equation}\label{formula fujiki}\int_{X}\alpha^{2n}=c_{X}B_{X}(\alpha,\alpha)^{n}\end{equation}
and for $0\neq\omega\in H^{0}(\Omega^{2}_{X})$
\begin{equation}B_{X}(\omega+\overline{\omega},\omega+\overline{\omega})>0.
\end{equation}\end{thm}
\begin{defn}\label{definizione forma beauville}
The quadratic form $B_{X}$ of the previous theorem is the
Beuville form.\\
The constant $c_{X}$ is the Fujiki constant.
\end{defn}
\begin{rem}The formula (\ref{formula fujiki}) is named the Fujiki formula;
its  polarized form is the following
\begin{equation}\label{formula fujiki polariizzata}
\int_{X}\alpha_{1}\wedge...\wedge\alpha_{2n}=\frac{c_{X}}{2n!}
\sum_{\sigma\in
S_{2n}}B(\alpha_{\sigma(1)},\alpha_{\sigma(2)})...B(\alpha_{\sigma(2n-1)},\alpha_{\sigma(2n)})
\end{equation}
\end{rem}
This section is organized as follows. In subsection 1 we simply
recall the basis of $H^{2}(\widetilde{\mathcal{M}},\mathbb{Q})$
given by O'Grady. Subsections 2, 3 and 4 are devoted to extract
from O'Grady's basis a basis of
$H^{2}(\widetilde{\mathcal{M}},\mathbb{Z})$. In subsection 2 it is
shown that the submodule
$\widetilde{\mu}(H^{2}(\mathcal{J},\mathbb{Z}))$ obtained  by
means of the Donaldson morphism (see Definition
\ref{definizionemutilde}) is saturated (Proposition \ref{sat}). In
subsection 3 it is proved that $2$ divides
$c_{1}(\widetilde{\Sigma})$ in
$H^{2}(\widetilde{\mathcal{M}},\mathbb{Z})$ (Theorem
\ref{divisibilita'}). Using these results and letting
$\widetilde{B}$ be the strict transform in
$\widetilde{\mathcal{M}}$ of the locus parametrizing stable
sheaves non locally free, we can easily prove  Theorem \ref{base
intera} of subsection 4, asserting that
\[H^{2}(\widetilde{\mathcal{M}},\mathbb{Z})=
\mathbb{Z}\frac{c_{1}(\widetilde{\Sigma})}{2}\oplus
\mathbb{Z}c_{1}(\widetilde{B})\oplus
\widetilde{\mu}(H^{2}(\mathcal{J},\mathbb{Z})).\] Finally  in
subsection 5 few intersection numbers on $\widetilde{\mathcal{M}}$
are needed to completely determine the Beauville form.
\subsection{The rational basis}
In  this subsection we simply recall the  basis of
$H^{2}(\widetilde{\mathcal{M}},\mathbb{Q})$ given by O'Grady.\\
There exists a smaller compactification, namely the Uhlenbeck
compactification, $\mathbf{M}^{U}$ of the $\mu-$stable locus of
$\mathbf{M}^{0}_{(2,0,-2)}$. Moreover $\mathbf{M}^{U}$ is endowed
with a surjective map
\[\varphi: \mathbf{M}^{0}_{(2,0,-2)}\longrightarrow \mathbf{M}^{U}. \]
Associated with $\varphi$ there is a cohomological linear map, the
Donaldson morphism
\begin{equation}\label{donaldson}\mu:H^{2}(\mathcal{J},\mathbb{Z}) \longrightarrow
H^{2}(\mathbf{M}^{U},\mathbb{Z})\end{equation} having the
following property (see \cite{Li93}, \cite{FM94} and \cite{Mo93}).\\
 Given $\mathcal{F}$ a flat family of sheaves in $\mathbf{M}^{0}_{(2,0,-2)}$,
 parametrized by a scheme $X$, letting
\[f_{\mathcal{F}}:X\rightarrow \mathbf{M}^{0}_{(2,0,-2)}\] be the modular
morphism and letting $p$ and $q$ be the two projections  to
$\mathcal{J}$ and $X$ respectively
\begin{equation}\label{propdonalson}
(f_{\mathcal{F}}^{*}\circ\varphi^{*}\circ \mu)(\alpha)= q_{*}(
p^{*}(\alpha)\cup c_{2}(\mathcal{F}))\end{equation}
for any $\alpha\in H^{2}(\mathcal{J},\mathbb{Z})$.\\
 We can finally define
\begin{defn}\label{definizionemutilde}
\[\widetilde{\mu}:=
\widetilde{\pi}^{*}\circ\varphi^{*}\circ\mu:H^{2}(\mathcal{J},\mathbb{Z})\rightarrow
H^{2}(\widetilde{\mathcal{M}},\mathbb{Z}).\]
\end{defn}
In order to recall the basis of
$H^{2}(\widetilde{\mathcal{M}},\mathbb{Q})$ we need to recall the
definition of a divisor on $\widetilde{\mathcal{M}}.$
\begin{defn} Let $B\subset \mathbf{M}^{0}_{(2,0,-2)}$
be the locally closed subset parametrizing non locally free stable
sheaves, O'Grady defines
\[\widetilde{B}:=\overline{\widetilde{\pi}^{-1}(B)}.\]
\end{defn}
We are now ready to present the rational basis: the following is
Proposition (7.3.3) of \cite{OG03}.
\begin{prop}
The homomorphism
$\widetilde{\mu}:H^{2}(\mathcal{J},\mathbb{Z})\rightarrow
H^{2}(\widetilde{\mathcal{M}},\mathbb{Z})$ is injective and \[
\begin{array}{ccc}
\widetilde{\mu}(H^{2}(\mathcal{J},\mathbb{Q})), &
\mathbb{Q}c_{1}(\widetilde{\Sigma}), & \mathbb{Q}c_{1}(\widetilde{B})\\
\end{array}\]
are linearly independent subspaces of
$H^{2}(\widetilde{\mathcal{M}},\mathbb{Q})$. Moreover
\[\widetilde{\mu}(H^{2}(\mathcal{J},\mathbb{Q}))\oplus
\mathbb{Q}c_{1}(\widetilde{\Sigma})\oplus
\mathbb{Q}c_{1}(\widetilde{B})=H^{2}(\widetilde{\mathcal{M}},\mathbb{Q})
\]
\end{prop}
 We will start
from the given basis to find $8$ independent generators of
$H^{2}(\widetilde{\mathcal{M}},\mathbb{Z})$.
\subsection{The image of the Donaldson's morphism}
Our first result  in the study of the 2-cohomology of
$\widetilde{\mathcal{M}}$ is the  following.
\begin{prop}\label{sat}The image of $\widetilde{\mu}:=
\widetilde{\pi}^{*}\circ\varphi^{*}\circ\mu:H^{2}(\mathcal{J},\mathbb{Z})
\longrightarrow H^{2}(\widetilde{\mathcal{M}},\mathbb{Z})$
 is a saturated submodule.
\end{prop}
\begin{proof}
By simple linear algebra it is enough to prove that there exists
\[\alpha:H^{2}(\widetilde{\mathcal{M}},\mathbb{Z})\longrightarrow\mathbb{Z}^{n}\]
such that the restriction of $\alpha$ to
$\widetilde{\mu}(H^{2}(\mathcal{J},\mathbb{Z}))$ is injective and
$\alpha(\widetilde{\mu}(H^{2}(\mathcal{J},\mathbb{Z})))$ is a
saturated submodule. Our map $\alpha$ will be the pull back map
associated to a  closed embedding $\widetilde{B}_{p}
\subset\widetilde{\mathcal{M}}$
 that we are going to define.\\
Given $p\in \mathcal{J}\setminus\mathcal{J}[2]$ let $B_{p}$ be the
locus of $\mathbf{M}_{(2,0,-2)}^{0}$ parametrizing non locally
free stable
sheaves, with singularities in $p$ or $-p$: $B_{p}$ is locally closed but not closed.
 We define:
\[\widetilde{B}_{p}:=\overline{\widetilde{\pi}^{-1}B_{p}}\]
namely $\widetilde{B}_{p}$ is the strict transform of $B_{p}$ in
$\widetilde{\mathcal{M}}$. $\widetilde{B}_{p}$ has already been
described in \cite{OG03}: the following proposition can be,
almost completely, extracted from  subsection 5.1 of that paper:
\begin{prop}\label{51og}\begin{enumerate}
\item{$B_{p}$ parametrizes
semistable sheaves $N$ fitting into exact sequences of the form
\begin{equation}\label{ddes}0\rightarrow N
\rightarrow
V\rightarrow\mathbb{C}_{p}\oplus\mathbb{C}_{-p}\longrightarrow
0\end{equation} where $V$ belongs to the moduli space $Muk$ of
Mukai-stable (see \cite{BDL01}) rank two vector bundles with
trivial determinant and $c_{2}^{hom}$, namely either $V=L\oplus
L^{\vee}$ with
$L\in\widehat{\mathcal{J}}\setminus\widehat{\mathcal{J}}[2]$ or
$V$ represents a non trivial element in $Ext^{1}(L,L)$ with $L\in
\widehat{\mathcal{J}}[2]$.} \item{The rational map\[\begin{matrix}
\zeta: & Muk &\rightarrow &
K(\widehat{\mathcal{J}})\\
& L\oplus L^{\vee} & \mapsto & (L, L^{\vee})\\
\end{matrix}\] extends to an isomorphism (here,
 $K(\widehat{\mathcal{J}})$ is the Kummer surface of
 $\widehat{\mathcal{J}}$, namely the locus of $Hilb^{2}(\widehat{\mathcal{J}})$ whose points
 correspond to schemes with associated cycles summing up to $0\in \widehat{\mathcal{J}}$ ).}
\item{ The rational map \[\phi: \widetilde{B}_{p}\rightarrow Muk\]
defined on  $\widetilde{\pi}^{-1}(B_{p})$  associating to
$x\in\widetilde{\pi}^{-1}([F])$ the class $[F^{\vee \vee}]$ in
$Muk$ extends to a regular morphism.} \item{The composition of the
extensions
\[\psi:= \zeta\circ\phi:\widetilde{B}_{p}\rightarrow K(\widehat{\mathcal{J}})\]
endows $\widetilde{B}_{p}$ with the structure of a
$\mathbb{P}^{1}$-bundle.} \item{Denote $\tilde{i}:
\widetilde{B}_{p} \rightarrow \widetilde{\mathcal{M}}$  the closed
embedding, then there exists a map $i_{U}$ making commutative the
following
 diagram
\[
\CD \widetilde{B}_{p}  @>\tilde{i}>> \widetilde{\mathcal{M}}\\
@V\psi VV @V\varphi\circ\widetilde{\pi} VV \\
K(\widehat{\mathcal{J}})  @>i_{U}>> \mathbf{M}^{U}.\\
\endCD
\]}

\end{enumerate}
\end{prop}
\begin{proof}
(1),(3) and (4) are proved in Lemma (4.3.3) and in section 5.1 of \cite{OG03}.\\
(2) is contained in  Theorem 5.6 of \cite{BDL01}.\\
(5) follows from the classification of the fibers of $\varphi$
(see \cite{HL97} ) since the fibers of $\psi$ are contracted by
$\varphi\circ\widetilde{\pi}$.
\end{proof}
The proof Proposition \ref{sat} will follow easily from the
following claim.
\begin{claim}\label{claimsat}Let
$b:Bl_{\mathcal{J}[2]}\mathcal{J}\rightarrow\mathcal{J}$ be the
blow up of the $2-torsion$ points of $\mathcal{J}$, let
$q:Bl_{\mathcal{J}[2]}\mathcal{J}\rightarrow K(\mathcal{J})$ be
the quotient by the involution, let
$e:K(\widehat{\mathcal{J}})\rightarrow K(\mathcal{J})$ be the
identification induced by the one between $\mathcal{J}$ and
$\widehat{\mathcal{J}}$, then
\[\widetilde{i}^{*}\circ\widetilde{\mu}(H^{2}(\mathcal{J},\mathbb{Z}))=
\psi^{*}\circ e^{*}\circ q_{*}\circ
b^{*}(H^{2}(\mathcal{J},\mathbb{Z}))\]
\end{claim}
\begin{proof}
The first step in proving the claim is to produce a 'complete'
family $ \mathcal{N}$ of
sheaves in $\widetilde{\pi}(\widetilde{B}_{p})$ and to study the topology of its base.\\
We start constructing a universal family $\mathcal{V}$
parametrized by $ K(\mathcal{J})$ for sheaves which are middle
terms of the sequences (\ref{ddes}): consider in fact the
structural sheaf of the tautological subvariety of $
Hilb^{2}(\mathcal{J})\times\mathcal{J}$ and restrict it to $
 K(\mathcal{J})\times\mathcal{J}$. This sheaf $\mathcal{G}$ can be seen as a
family of sheaves of length 2 quotients of
$\mathcal{O}_{\mathcal{J}}$, they fit in exact sequences of the
form
\[0\rightarrow\mathbb{C}_{x}\rightarrow G\rightarrow\mathbb{C}_{-x}\rightarrow0\]
( $x$ is a point of $\mathcal{J}$), moreover such a sequence
splits if and only if $x\neq -x$. Recall then that for any
$x\in\mathcal{J}$ the sheaf $\mathbb{C}_{x}$ satisfies the weak
index theorem \cite{Mu81} , therefore any $G$,  being the middle
term of such an extension, satisfies W.I.T. By the theory of
\cite{Mu87}, applying the Fourier-Mukai transform to the family
$\mathcal{G}$, we obtain a new family $\mathcal{V}$  of sheaves on
$\mathcal{J}$. Since the Fourier-Mukai transform is a
self-equivalence of the derived category of $\mathcal{J}$ and the
Fourier-Mukai transform of a sheaf $\mathbb{C}_{x}$ is
$\mathcal{O}(\Theta_{x}-\Theta)$, $\mathcal{V}$  is a family
parametrizing
bijectively sheaves in $Muk$.\\
 We can now construct the family $\mathcal{N}$: consider the relative
Quot-scheme parametrizing 0-dimensional length 2 quotients of the
`fibers' of $\mathcal{V}$, take the closed subscheme $Q$ of the
sheaves whose support is exactly \{$p,-p$\} and, in this
subscheme, the open $\mathcal{U}$ corresponding to quotients
having semistable kernels. On $\mathcal{J}\times\mathcal{U}$ we
have the following exact sequence of $\mathcal{U}$-flat sheaves:
\begin{equation}\label{seqkft}0\longrightarrow\mathcal{N}\longrightarrow (id\times pr)
^{*}\mathcal{V}\longrightarrow\mathcal{T}\longrightarrow
0\end{equation} where $pr$ is the natural projection from
$\mathcal{U}$ to $ K(\mathcal{J})$ and $\mathcal{T}$ and
$\mathcal{N}$ respectively the restrictions of the tautological
quotient and the tautological kernel of the previous Quot-scheme.
More explicitly $Q$ can be obtained as follows. Let $q_{i}$ be the
projection of $\mathcal{J}\times  K(\mathcal{J})$ to its $i$-th
factor: then
$q_{2*}\mathcal{H}om(\mathcal{V},q_{1}^{*}\mathbb{C}_{p})$ and
$q_{2*}\mathcal{H}om(\mathcal{V},q_{1}^{*}\mathbb{C}_{-p})$ are
rank 2-vector bundles on $ K(\mathcal{J})$ and  there is an
identification
\[Q=\mathbb{P}(q_{2*}\mathcal{H}om(\mathcal{V},q_{1}^{*}\mathbb{C}_{p}))
\times_{K(\mathcal{J})}
\mathbb{P}(q_{2*}\mathcal{H}om(\mathcal{V},q_{1}^{*}\mathbb{C}_{-p})).\]
In particular the fiber on $V\in Muk\simeq K(\mathcal{J})$ is
$\mathbb{P}(Hom(V,\mathbb{C}_{p}))\times\mathbb{P}(Hom(V,\mathbb{C}_{-p}))\simeq
\mathbb{P}^{1}\times \mathbb{P}^{1}$.
 To understand $\mathcal{U}$ we have to detect
the unstable locus of each fiber. Let $([l_{p}],[l_{-p}])\in
\mathbb{P}(Hom(V,\mathbb{C}_{p}))\times\mathbb{P}(Hom(V,\mathbb{C}_{-p}))$,
then $([l_{p}],[l_{-p}])$ gives an unstable $N$ if and only if
there exists $L\subset V$ such that $c_{1}(L)=0$ with
$l_{p}(L_{p})=l_{-p}(L_{-p})=0$. It can be easily checked that, if
$V$ is a direct sum there are only 2 points giving $N$ unstable
and if $V$ is a non trivial extension there is a unique point
giving $N$ unstable: therefore $\mathcal{U}$ is  isomorphic
to the complement of a 2 dimensional variety in the 4-fold $Q$.\\
The second step of the proof of the claim consists in computing
the pull-back of
$\varphi^{*}\circ\mu(H^{2}(\mathcal{J},\mathbb{Z}))$  via the
modular map $f_{\mathcal{N}}$ associated to the family
$\mathcal{N}$. The result is the following:
\begin{equation}\label{semiclaim}f_{\mathcal{N}}^{*}\circ\varphi^{*}\circ \mu
(H^{2}(\mathcal{J},\mathbb{Z})) = pr^{*}\circ q_{*} \circ b^{*}
(H^{2}(\mathcal{J},\mathbb{Z})).
\end{equation}
Let $p_{i}$  be the projection from $\mathcal{J}\times\mathcal{U}$
to the i-th factor, and apply Whitney's formula to the exact
sequence (\ref{seqkft}), since the support of $\mathcal{T}$ is
$p_{1}^{-1}(p)\cup p_{1}^{-1}(-p)$, hence $c_{1}(\mathcal{T})=0$,
we get:
\[c_{2}  (\mathcal{N})=c_{2}((id\times pr)^{*}\mathcal{V})-c_{2}
(\mathcal{T}).\]
 So for any
$\alpha$ belonging to $H^{2}(\mathcal{J},\mathbb{Z})$, we have
\[p_{2*}(p_{1}^{*}{\alpha}\cup c_{2}(\mathcal{N}))=
 p_{2*}(p_{1}^{*}{\alpha}\cup
 (c_{2}((id\times pr)^{*}\mathcal{V})-c_{2}(\mathcal{T})))=
 p_{2*}(p_{1}^{*}{\alpha}\cup
 c_{2} ((id\times pr)^{*}\mathcal{V}))\]
 where the second equality is verified because
 $c_{2}(\mathcal{T})$ is the pull back of a 4-form
 from $\mathcal{J}$. Since, obviously, $p_{i}=q_{i}\circ (id\times pr)$
we can simplify the last term of the last equation:
\[p_{2*}(p_{1}^{*}{\alpha}\cup
 c_{2} ((id\times pr)^{*}\mathcal{V}))=
 p_{2*}((id\times pr)^{*}(q_{1}^{*}\alpha\cup
 c_{2}(\mathcal{V})))=
 pr^{*}q_{2*}(q_{1}^{*}\alpha\cup c_{2}(\mathcal{V})).\]
 Call now $\mathcal{I}$ the incidence subvariety of
$\mathcal{J}\times K(\mathcal{J})$: as we said earlier,
$\mathcal{V}$ is the Fourier-Mukai transform of the family
$\mathcal{O}_{\mathcal{I}}$; let then $FM^{\bullet}$ be  the
cohomological Fourier-Mukai transform on
$H^{\bullet}(\mathcal{J},\mathbb{Z})$, it also acts on
$H^{\bullet}(\mathcal{J}\times K(\mathcal{J}),\mathbb{Z})$ by
means of  the K\"unneth decomposition. Since, by
Grothendieck-Riemann-Roch,
$FM^{\bullet}(ch(\mathcal{O}_{\mathcal{I}}))=ch(\mathcal{V})$, and
moreover $FM^{\bullet}$ acts on $H^{2}(\mathcal{J},\mathbb{Z})$ as
an isometry $FM$( with respect to the intersection form ( see
\cite{Mu87} )) we get
\[ q_{2*}(q_{1}^{*}\alpha\cup c_{2}(\mathcal{V}))=-[q_{2*}(q_{1}^{*}\alpha\cup ch(\mathcal{V}))]_{2}=
-[ q_{2*}(q_{1}^{*}\circ FM(\alpha)\cup
 ch(\mathcal{O}_{\mathcal{I}}))]_{2}\] where the first equality is true since $\mathcal{V}$
parametrizes sheaves with trivial determinant.\\
But by Grothendieck-Riemann-Roch
$ch_{2}(\mathcal{O}_{\mathcal{I}})$ is represented exactly by the
incidence variety which can be identified, in an obvious way, to
$Bl_{\mathcal{J}[2]}\mathcal{J}$:  the restrictions of the two
projections become then  the maps $b$ and $q$ of the claim and
then we obtain
\[[q_{2*}(q_{1}^{*}\circ FM(\alpha)\cup
ch(\mathcal{O}_{\mathcal{I}}))]_{2} =q_{*}(b^{*}\circ
FM(\alpha)),\] consequently replacing this in the previous
equations we get the formula
\[p_{2*}(p_{1}^{*}{\alpha}\cup c_{2}(\mathcal{N}))=-pr^{*}(q_{*}(b^{*}\circ FM(\alpha)))\]
and finally by the property (\ref{propdonalson}) of the
Donaldson's morphism
\[f_{\mathcal{N}}^{*}\circ\varphi^{*}\circ\mu(\alpha)=-pr^{*}(q_{*}(b^{*}\circ
FM(\alpha)))\] which implies (\ref{semiclaim}) since $FM$ is in
particular an integral isomorphism on
$H^{2}(\mathcal{J},\mathbb{Z})$.\\
To complete the proof of the claim it remains to relate  the map
$\psi:\widetilde{B}_{p}\rightarrow K(\widehat{\mathcal{J}})$ of
Proposition \ref{51og} (see
the diagram) to $pr:\mathcal{U}\rightarrow K(\mathcal{J})$.\\
Notice that, by Theorem 8.2.11 of \cite{HL97}, $\varphi\circ
f_{N}$ is constant on the fibers of $pr$ and, since $\mathcal{U}$
is an open in a locally trivial bundle on a smooth base, it
factors through the base. More precisely, using the identification
$ e:K(\widehat{\mathcal{J}})\rightarrow K(\mathcal{J}) $ induced
by $\mathcal{J}\simeq\widehat{\mathcal{J}}$ it can be easily seen
that:
\[\varphi\circ f_{\mathcal{N}}=i_{U}\circ e^{-1} \circ pr.\] Therefore by
(\ref{semiclaim})
\[pr^{*}\circ q_{*} \circ b^{*}
(H^{2}(\mathcal{J},\mathbb{Z}))=pr^{*}\circ (e^{-1})^{*} \circ
i_{U}^{*} \circ \mu (H^{2}(\mathcal{J},\mathbb{Z})).\] Since the
complement of $\mathcal{U}$ has complex codimension 2 in the
$\mathbb{P}^{1}\times\mathbb{P}^{1}$ bundle $Q$,
$pr^{*}:H^{2}(K(\mathcal{J}),\mathbb{Z})\rightarrow
H^{2}(\mathcal{U},\mathbb{Z})$ is injective, therefore $i_{U}^{*}
\circ\mu(H^{2}(\mathcal{J},\mathbb{Z})) =e^{*}\circ q_{*} \circ
b^{*} (H^{2}(\mathcal{J},\mathbb{Z}))$ and by the commutativity of
the diagram in (5) of \ref{51og}\[\widetilde{i}^{*}\circ
\widetilde{\mu}(H^{2}(\mathcal{J},\mathbb{Z}))=\psi^{*}\circ
i_{U}^{*} \circ \mu (H^{2}(\mathcal{J},\mathbb{Z}))=\psi^{*}\circ
e^{*}\circ q_{*} \circ b^{*} (H^{2}(\mathcal{J},\mathbb{Z})).\]
\end{proof}

Now we finish proving Proposition \ref{sat}: it is well known that
$q_{*}\circ b^{*}(H^{2}(\mathcal{J},\mathbb{Z}))$ is the
orthogonal, with respect to the intersection form, of the
submodule generated by the nodal classes in $H^{2}(
K(\mathcal{J}),\mathbb{Z})$ \cite{BPV84}, hence it is saturated;
since $\psi:\widetilde{B}_{p}\rightarrow K(\widehat{\mathcal{J}})$
is a $\mathbb{P}^{1}$-bundle (see (4) of Proposition \ref{51og}),
$\psi^{*}\circ e^{*}\circ q_{*}\circ
b^{*}(H^{2}(\mathcal{J},\mathbb{Z}))$ is saturated too, thus the
claim implies the proposition.\end{proof}

\subsection{The divisibility of $\widetilde{\Sigma}$ }

The goal of this subsection is to prove the following theorem
\begin{thm}\label{divisibilita'}
There exists $A\in H^{2}(\widetilde{\mathcal{M}},\mathbb{Z})$ such
that $2A=c_{1}(\widetilde{\Sigma})$.
\end{thm}
In Proposition \ref{v02-2} we  have exhibited a birational map
$b\colon\widetilde{\mathcal{M}}\rightarrow\widetilde{\mathbf{M}}_{(0,2\Theta,2)}^{0}$
furthermore we have proved  that
$b^{*}(\widetilde{\Sigma}_{(0,2\Theta,-2)}^{0})=\widetilde{\Sigma}$:
since $\widetilde{\mathcal{M}}$ and
$\widetilde{\mathbf{M}}_{(0,2\Theta,2)}^{0}$ are both symplectic,
they are isomorphic (via $b$) in codimension 1. It follows that
$b$ induces isomorphisms on Picard groups and integral
2-cohomology groups. Theorem \ref{divisibilita'} is therefore a
consequence  of the following proposition.
\begin{prop}\label{divis 02-2}
$2|c_{1}(\widetilde{\Sigma}^{0}_{(0,2\Theta,-2)})$ in
$H^{2}(\widetilde{\mathbf{M}}_{(0,2\Theta,-2)}^{0},\mathbb{Z})$.
\end{prop}
\noindent{\em Idea of the Proof.} We will prove the desired
divisibility of $c_{1}(\Sigma_{(0,2\Theta,-2)}^{0})$ using
particular features of
$\widetilde{\mathbf{M}}_{(0,2\Theta,-2)}^{0}$.
$\widetilde{\mathbf{M}}_{(0,2\Theta,-2)}^{0}$  is the symplectic
desingularization of a moduli space
$\mathbf{M}_{(0,2\Theta,-2)}^{0}$ whose general point parametrizes
a sheaf of the form $F=(j_{D})_{*}L$, where $j_{D}: D\rightarrow
\mathcal{J}$ and $L$ is a line bundle on $D$ such that
$-id^{*}(F)=F$ (recall that $-id^{*}$ fixes every $D\in
|2\Theta|$) (see \ref{-1invarianza}). If the support $D$ of $F$
does not pass through a 2-torsion point of $\mathcal{J}$, $F$ is
the pull back of a sheaf on a curve on the Kummer surface
$Kum:=K(\mathcal{J})$ associated to $\mathcal{J}$. This  shows the
existence of a dominant rational map $\tau$  from a moduli space
of sheaves on $Kum$, ($\mathbf{M}_{(0,d^{*}H,-1)}(Kum)$ defined in
\ref{def mkum e ukum}) to
$\mathbf{M}_{(0,2\Theta,-2)}^{0}$.\\
The map $\tau$ is studied in Proposition \ref{2:1} and is proved
to be generically [2:1].\\
The idea of the proof of the divisibility of
$c_{1}(\widetilde{\Sigma}_{(0,2\Theta,-2)}^{0})$  is  to extend
$\tau$ to a finite map on a big (namely having complement of
codimension strictly bigger than 1) open subset of
$\widetilde{\mathbf{M}}_{(0,2\Theta,-2)}^{0}$ and then  study (the
closure of) its branch locus $R$: the general theory of
 double coverings implies in fact that $c_{1}(R)$ is
2-divisible in
$H^{2}(\widetilde{\mathbf{M}}_{(0,2\Theta,-2)}^{0},\mathbb{Z})$.
In order to control the branch locus of the extension of $\tau$
(exhibited in Proposition \ref{completamento2a1}) we  need to
study the complement of the image in
$\widetilde{\mathbf{M}}_{(0,2\Theta,-2)}^{0}$ of the original map
$\tau$ (see Proposition \ref{mult}) and we need to recall some
basic facts about the action of $\mathcal{J}[2]$ on $\mathcal{J}$ and
$Kum$ (see Remark \ref{rem j[2]}). The branch locus $R$ is
studied, not determined, in the proof of the divisibility of
$c_{1}(\widetilde{\Sigma}_{(0,2\Theta,-2)}^{0})$ : it surely
satisfies $R=\widetilde{\Sigma}_{(0,2\Theta,-2)}^{0}+ Q$, $Q$
being a divisor such that $2|c_{1}(Q)\in
H^{2}(\widetilde{\mathbf{M}}_{(0,2\Theta,-2)}^{0},\mathbb{Z})$.
This implies the divisibility of
$c_{1}(\widetilde{\Sigma}_{(0,2\Theta,-2)}^{0})$.\qed\medskip\\
Thus we are reduced to prove the 2-divisibility of
$c_{1}(\Sigma_{(0,2\Theta,-2)}^{0})$ in the group
$H^{2}(\widetilde{\mathbf{M}}_{(0,2\Theta,-2)}^{0},\mathbb{Z})$:
the advantage of $\widetilde{\mathbf{M}}_{(0,2\Theta,-2)}^{0}$ is
that it can be easily related to a moduli space of sheaves on
\[Kum:=K(\mathcal{J}).\] To explain this relation we fix
notation and recall some well known classical result.
\begin{notation}
Let  $f_{|2\Theta|}:\mathcal{J}\rightarrow |2\Theta|^{\vee}$ be
the regular map associated with the complete linear system
$|2\Theta|$.\\  As in the proof of Lemma \ref{classificazione
curvelem} denote by $Kum_{s}:=f_{|2\Theta|}(\mathcal{J})$ the
singular Kummer surface associated with $\mathcal{J}$ and let
$f:\mathcal{J}\rightarrow Kum_{s}$ be the map induced by
$f_{|2\Theta|}$. It is well known that $f$ is identified with the
quotient of $\mathcal{J}$ by the involution $-id$, in particular
$f$ sends the $16$ $2$-torsion points of $\mathcal{J}$ to the 16
singular points of
$Kum_{s}$ and outside of them is an unramified double covering.\\
Let $d: Kum\rightarrow Kum_{s}$ be the desingularization map, and
$E:=\sum_{i=1}^{16}$ be the exceptional divisor of $d$. The
$E_{i}$'s are smooth rational curves on $Kum$, also known as nodal
curves.\end{notation}
\begin{defn}\label{def mkum e ukum} \begin{enumerate}\item{
Let $H$ be a plane in $ |2\Theta|^{\vee}$ and let
$D:=d^{*}H-\epsilon\sum_{i}E_{i}$ be an ample divisor in
$Pic(Kum)$, we denote by $\mathbf{M}_{(0,d^{*}H,-1)}(Kum)$  the
moduli space of sheaves on $Kum$ having Mukai's vector
$(0,d^{*}H,-1)$ and being semistable with respect to the
polarization $D$.}\item{We denote by
$U_{(0,d^{*}H,-1)}(Kum)\subset\mathbf{M}_{(0,d^{*}H,-1)}(Kum)$ the
open subscheme parametrizing sheaves whose Fitting subschemes (see
Remark \ref{fitting astra}) do not intersect the nodal curves.}
\end{enumerate}\end{defn}
\begin{defn}\label{def V e Vs}\begin{enumerate}\item{We denote
by $V_{(0,2\Theta,-2)}:=\Phi^{-1}(S\cup
R(1))\subset\mathbf{M}_{(0,2\Theta,-2)}^{0}$ the open subscheme
parametrizing sheaves whose Fitting subscheme does not pass
through 2-torsion points.} \item{We denote by $V^{s}\subset
V_{(0,2\Theta,-2)}$ the open subscheme parametrizing stable
sheaves.}
\end{enumerate}\end{defn}
 We can now begin to study the relation between
$\mathbf{M}_{(0,d^{*}H,-1)}(Kum)$ and
$\mathbf{M}_{(0,2\Theta,-2)}^{0}$
\begin{prop}\label{2:1}
The functor $f^{*}\circ d_{*}=b_{*}\circ q^{*}:\mathcal{C}oh(
Kum)\rightarrow\mathcal{C}oh(\mathcal{J})$ ($b$ and $q$ as in
claim \ref{claimsat}) induces a regular map
\[\begin{matrix}  \tau : & U_{(0,d^{*}H,-1)}(Kum)
&\rightarrow & \mathbf{M}_{(0,2\Theta,-2)}^{0}\\
&[F] & \mapsto & [(f^{*}\circ d_{*})F]
\end{matrix}\]such that \begin{enumerate}\item{$ \tau (
U_{(0,d^{*}H,-1)}(Kum))= V_{(0,2\Theta,-2)}$} \item{ If $x\in
V^{s}$ then $ \tau ^{-1}(x)$ consists of two distinct points.}
\item{If $x\in V_{(0,2\Theta,-2)}\setminus V^{s}$ then $ \tau
^{-1}(x)$ consists of a point.}\end{enumerate}
\end{prop}
\begin{proof}
Let $[F]\in U_{(0,d^{*}H,-1)}(Kum)$ and let $C_{K}\in|d^{*}(H)|$
be its Fitting subscheme, then $C:=f^{*}\circ d_{*}(C_{K})$
belongs to $\Phi^{-1}(S\cup R(1))$ (see Proposition
\ref{classificazione curve}). Since $C_{K}$ is the quotient of $C$
by the restriction of $-1$, it is easily verified that $b_{*}\circ
q^{*}(F)\in \mathbf{M}_{(0,2\Theta,-2\eta)}^{0}$ if
and only if $b_{*}\circ q^{*}(F)$ is semistable.\\
If $C\in S$ then the restriction of $F$ to $C_{K}$ is a line
bundle and the same property holds for $f^{*}\circ d_{*}(F)$
making it a stable sheaf. If $C\in R(1)$ and the restriction of
$F$ to $C_{K}$ is a line bundle (of degree 1), then $f^{*}\circ
d_{*}(F)$ is a line bundle having  the same degree (1) on each
component of $C$: also in this case $f^{*}\circ d_{*}(F)$ can be
proved to be stable. Finally if $C\in R(1)$ and the restriction of
$F$ to $C_{K}$ is not a line bundle then $F$ is the push-forward
of a degree-0 line bundle $L$ from the normalization (isomorphic
to $C^{0}$) of $C_{K}$, it follows that $f^{*}\circ d_{*}(F)$ is
the push-forward,  from the desingularization of $C$, of a sheaf
restricting to $L$ on each connected component: it is therefore
polystable. Thus $\tau$ is regular on $U_{(0,d^{*}H,-1)}(Kum)$.\\
Furthermore, for $C\in S\cup R(1)$, $\Phi^{-1}(C)$ is irreducible.
In the first case, as explained in the proof of Proposition
\ref{caratteristica fibre}, $\Phi^{-1}(C)$ is identified with the
kernel of the map $a:J(C)\rightarrow \mathcal{J}$ induced by the
embedding $i:C\rightarrow \mathcal{J}$. To show that $Ker(a)$ is
an irreducible torus it is enough to check that
$i_{*}:H_{1}(C,\mathbb{Z})\rightarrow
H_{1}(\mathcal{J},\mathbb{Z})$ is surjective: this follows from
Lefschetz hyperplane theorem since $C\in |2\Theta|$ and $\Theta $
is
ample.\\
In the second case, the stable locus of $\Phi^{-1}(C)$ is
isomorphic to a $\mathbb{C}^{*}$ bundle over $\mathcal{J}$ (see
the proof of Proposition \ref{caratteristica fibre}) and it is
easily seen to be dense in $\Phi^{-1}(C)$.\\
For $C\in S\cup R(1)$, let $M_{C_{K}}\subset
U_{(0,d^{*}H,-1)}(Kum)$ be the locus parametrizing sheaves having
$C_{K}$ as Fitting subscheme. The restriction of $\tau$ induces a
map $\tau_{C}\colon M_{C_{K}}\rightarrow \Phi^{-1}(C)$. $\tau_{C}$
is generically [2:1], in fact if $F_{1}$ and $F_{2}$ are sheaves
restricting to degree-1 line bundle on $C_{K}$, then $f^{*}\circ
d_{*}(F_{1})\simeq f^{*}\circ d_{*}(F_{2})$ if and only if they
differ by tensor product with $(i_{K})_{*}L$, where
$i_{K}:C_{K}\rightarrow Kum$ is the closed embedding and $L$ is
the degree 0 line bundle on $C_{K}$ associated with the unramified
double covering $C$.\\
Since $M_{C_{K}}$ and $\Phi^{-1}(C)$ are projective and have the
same dimension  $\tau_{C}$ is surjective: this implies item 1.\\
Item 2 follows since $V^{s}$ parametrizes sheaves having locally
free restriction to their support (see the proof of Proposition
\ref{caratteristica fibre}), hence they must be images via $\tau$
of sheaves having the same property.\\
Item 3 can be checked directly.

\end{proof}
\begin{rem}\label{-1invarianza}
Proposition \ref{2:1} shows that the general sheaf $F$ such that
$[F]\in\mathbf{M}^{0}_{(0,2\Theta,-2)}$ is the pull-back of a
sheaf on $Kum_{s}$, in particular $(-id)^{*}(F)=F$, since
$\mathbf{M}^{0}_{(0,2\Theta,-2)}$ is irreducible,
$\mathbf{M}^{0}_{(0,2\Theta,-2)}$ is in the fixed locus of the map
induced on $\mathbf{M}_{(0,2\Theta,-2)}$ by $(-id)^{*}$

\end{rem}
\begin{defn}\label{defuskum}We define $U^{s}(Kum):=  \tau  ^{-1}(V^{s})\subset
U_{(0,d^{*}H,-1)}(Kum)$ and denote by
\[  \tau  ^{s}:U^{s}(Kum)\rightarrow V^{s}\] the restriction of $  \tau  $.
\end{defn}
\begin{rem}The proof of Proposition \ref{2:1} shows that
$U^{s}(Kum)\subset U_{(0,d^{*}H,-1)}(Kum)$ is the locus
parametrizing sheaves having locally free restrictions to their
supports.
\end{rem}
The following corollary extracts from Proposition \ref{2:1} what
we will use to prove the divisibility of
$c_{1}(\widetilde{\Sigma}^{0}_{(0,2\Theta,-2)})$.
\begin{cor}\label{etale 2a1}
The map $  \tau  ^{s}:U^{s}(Kum)\rightarrow V^{s}$ is proper
\'etale [2:1].
\end{cor}
\begin{proof}Since $U^{s}(Kum)$ and $V^{s}$ parametrize stable sheaves
they are smooth and have pure dimension 6. Corollary \ref{etale
2a1} follows
  since, from Proposition \ref{2:1}, each fiber of
   $\tau ^{s}$  consists of 2 distinct points.
\end{proof}
Since the O'Grady desingularization doesn't modify the stable
locus, $V^{s}$ is naturally included in
$\widetilde{\mathbf{M}}_{(0,2\Theta,-2)}^{0}$:
 we now analyze  $( \widetilde{\mathbf{M}}_{(0,2\Theta,-2)}^{0}\setminus V^{s})$.
\begin{prop}\label{per irriducibilita'}
Set
$\widetilde{\Phi}:=\Phi\circ\widetilde{\pi}^{0}_{(0,2\Theta,-2)}\colon
\widetilde{\mathbf{M}}_{(0,2\Theta,-2)}^{0}\rightarrow |2\Theta|$.
For $x\in\mathcal{J}[2]$, let $H_{x}$ be the plane in $|2\Theta|$
corresponding to curves passing through $x$, then the locus
$\Delta_{x}:=\widetilde{\Phi}^{-1}H_{x}\subset\widetilde{\mathbf{M}}_{(0,2\Theta,-2)}^{0}$
is an irreducible divisor.
\end{prop}
This proposition is a consequence of the following lemma
\begin{lem}\label{irrid per c in N(1)}
If $C\in N(1)$ then $\Phi^{-1}(C)$ is irreducible.
\end{lem}
\begin{proof}
As explained in the proof of Proposition \ref{caratteristica
fibre} the open dense subset of $\Phi^{-1}(C)$,whose points
parametrize sheaves restricting to line bundles over $C$, is a
$\mathbb{C}^{*}$ bundle over the kernel of the map $a\colon
J(\widetilde{C})\rightarrow \mathcal{J}$ induced on Albanese
varieties by the map $h: \widetilde{C}\rightarrow \mathcal{J}$
obtained composing the normalization map of $C$ with its inclusion
in $\mathcal{J}$.\\The irreducibility of $\Phi^{-1}(C)$  is again
a consequence of the surjectivity of
$h_{*}:H_{1}(\widetilde{C},\mathbb{Z})\rightarrow
H_{1}(\mathcal{J},\mathbb{Z})$: since it can be checked that the
natural image of $\widetilde{C}$ in the blow up
$Bl_{x}\mathcal{J}$ of $\mathcal{J}$ in $x$ is an ample divisor
and moreover the blow up map induces an isomorphism on integral
1-homology groups, the desired surjectivity follows again by
Lefschetz hyperplane theorem.
\end{proof}
\noindent{\em Proof of Proposition \ref{per irriducibilita'}.}
Since $H_{x}\cap (S\cup R(1))=\emptyset$ we have,
\[\widetilde{\Phi}^{-1}H_{x}=\widetilde{\Phi}^{-1}(H_{x}\cap N(1))
\cup \widetilde{\Phi}^{-1}(H_{x}\cap (N(2)\cup N(3)\cup R(2)\cup
D)).\] But $\widetilde{\Phi}^{-1}(H_{x}\cap
N(1))=\Phi^{-1}(H_{x}\cap N(1))$ and by Lemma \ref{irrid per c in
N(1)} the latter is a fibration with irreducible 3 dimensional
fiber and irreducible 2 dimensional base: hence it is a
codimension 1 locally closed
subset.\\
On the other hand $\widetilde{\Phi}^{-1}(H_{x}\cap (N(2)\cup
N(3)\cup R(2)\cup D))$ is the inverse image of a codimension 2
subset via a Lagrangian fibration on an irreducible symplectic
variety: by Theorem 2 (step 5) of \cite{Ma99} it cannot contain a
codimension 1 subset. It follows that the divisor
$\Delta_{x}:=\widetilde{\Phi}^{-1}H_{x}$ is irreducible.\qed
\medskip\\ Finally we determine the multiplicity of
$\widetilde{\Phi}^{*}H_{x}$
\begin{prop}\label{mult}\[\widetilde{\Phi}^{*}H_{x}= \Delta_{x}.\]
\end{prop}
\begin{proof} Let $H\subset |2\Theta|$ be a
generic plane: we will show that
$c_{1}(\widetilde{\Phi}^{*}H)=c_{1}(\widetilde{\Phi}^{*}H_{x})$ is
indivisible in
$H^{2}(\widetilde{\mathbf{M}}_{(0,2\Theta,-2)}^{0},\mathbb{Z})$.
Since $\widetilde{\Phi}^{*}H$ is a reduced irreducible divisor  we
will prove the proposition exhibiting a curve intersecting it
transversally in a unique point. Consider the curves in
$|2\Theta|$ passing through 4 fixed points $\{a, -a, b, -b\}$
outside $\mathcal{J}[2]$: they form line $L$ in $|2\Theta|$, this
line consists of the inverse images of the hyperplane sections of
the singular Kummer surface passing through 2 fixed points. It is
easy to produce a flat family $\mathcal{F}$ on $\mathcal{J}\times
L$ such that, letting $j:C\rightarrow\mathcal{J}$ be the inclusion
of $C\in|2\Theta|$, for all $C\in L$
\[\mathcal{F}_{|\mathcal{J}\times
C}=j_{*}\mathcal{O}(j^{-1}(a)+j^{-1}(-a)).\] This family induces a
modular map $f_{\mathcal{F}}$ such that $\Phi\circ
f_{\mathcal{F}}$ is the identity map. For $H$ general, the image
of $f_{\mathcal{F}}$ obviously intersects $\Phi^{-1}H$
transversally in  a unique point: it follows that also the image
of the  lifting of $f_{\mathcal{F}}$ to
$\widetilde{\mathbf{M}}_{(0,2\Theta,-2)}^{0}$
intersects $\widetilde{\Phi}^{-1}H$ transversally in a unique
point.
\end{proof}
 We are now going to prove that $c_{1}(\widetilde{\Sigma}_{(0,2\Theta,-2)}^{0}) \in
H^{2}(\widetilde{\mathbf{M}}_{(0,2\Theta,-2)}^{0},\mathbb{Z})$ is
a class divisible by $2$: it will follow from the existence of a
square root of the class
$[\widetilde{\Sigma}_{(0,2\Theta,-2)}^{0}]$
in $Pic(\widetilde{\mathbf{M}}_{(0,2\Theta,-2)}^{0})$.\\
We will give a partial completion of the \'etale double covering
of Corollary \ref{etale 2a1}  and studying its branch locus we
will obtain the existence of such a square root: $ \tau  $ has to
be completed since its image does not contain the divisors
$\Delta_{x}$ and $\widetilde{\Sigma}^{0}_{(0,2\Theta,-2)}$ and
therefore it cannot give  relations on
 $Pic(\widetilde{\mathbf{M}}_{(0,2\Theta,-2)}^{0})$.\\
We anticipate that we will not be able to say whether our new double
covering ramifies or not on the $\Delta_{x}$'s, but it will be
enough to state that if it ramifies on a $\Delta_{x}$ it ramifies
on each $\Delta_{x}$: this will be achieved by means of a
$\mathcal{J}[2]$-action, on both the domain and the codomain, which
permutes the
$\Delta_{x}$ and is compatible with the double covering.\\
In the following remark we define and explain the actions that we
will use.
\begin{rem} \label{rem j[2]}\begin{enumerate}
\item{ Let $x$ be a 2-torsion point of $\mathcal{J}$,
then obviously
\[(-id)\circ t_{x}=t_{x}\circ(-id):\mathcal{J}\rightarrow\mathcal{J}\]
and this implies that the action by translation of
$\mathcal{J}[2]$ on $\mathcal{J}$ descends to an action on the
quotient $Kum_{s}$. Finally this action lifts to the blow up
$Kum$ of the singular locus of $Kum_{s}$.}
 \item{The $\mathcal{J}[2]$
action on $\mathcal{J}$ induces  $\mathcal{J}[2]$ actions on
$\mathbf{M}_{(0,2\Theta,-2)}^{0}$, on
$\Sigma_{(0,2\Theta,-2)}^{0}$ and on
$\Omega_{(0,2\Theta,-2)}^{0}$. In fact $\mathcal{J}[2]$ obviously
acts on $\mathbf{M}_{(0,2\Theta,-2)}$, on
$\Sigma_{(0,2\Theta,-2)}$, and on $\Omega_{(0,2\Theta,-2)}$;
furthermore, since for $x\in\mathcal{J}[2]$ $t_{x}^{*}2\Theta=
2\Theta$ and for $[F]\in\mathbf{M}_{(0,2\Theta,-2)}$ $\sum
c_{2}(t_{x}^{*}(F))=deg (c_{2}(F))x+\sum c_{2}(F)$ ,
$\mathbf{M}_{(0,2\Theta,-2)}^{0}$ is invariant for this action.\\
Moreover since $\mathcal{J}[2]$ permutes the 2-torsion points on
$\mathcal{J}$, the open subscheme $V_{(0,2\Theta,-2)}$ defined in
\ref{def V e Vs} is invariant under the $\mathcal{J}[2]$ action
and, since translations on $\mathcal{J}$ do not change stability,
the open subscheme $V^{s}\subset V_{(0,2\Theta,-2)}$ also defined
in \ref{def V e Vs} is invariant too.} \item{The $\mathcal{J}[2]$
action on $Kum$ induces an action on
$\mathbf{M}_{(0,d^{*}H,-1)}(Kum)$. In fact, since we choose a
polarization $D=d^{*}H-\epsilon\sum E_{i}$ (see \ref{def mkum e
ukum}) symmetric with respect to this action, the pull-backs by
automorphisms induced by $\mathcal{J}[2]$ do not change stability
and semistability. Moreover, since the $\mathcal{J}[2]$ action on
$Kum_{s}$ permutes the singular points, the open subscheme
$U_{(0,d^{*}H,-1)}(Kum)$ defined in \ref{def mkum e ukum} is
invariant under this action and, since the pull back of a sheaf,
locally free on its support, has the same property, the open
subscheme $U^{s}(Kum)\subset U_{(0,d^{*}H,1)}(Kum)$ defined in
\ref{defuskum} is invariant too .} \item{Since by 1) of this
remark the actions of $\mathcal{J}[2]$ commute with
$f:\mathcal{J}\rightarrow Kum_{s}$ and $d: Kum_{s}\rightarrow
Kum$, the pull backs of sheaves via the automorphisms induced by a
fixed $x\in\mathcal{J}[2]$  on $\mathcal{J}$ and $Kum$ commute
with $f_{*}\circ d^{*}$, therefore the regular morphism
\[  \tau  :U_{(0,d^{*}H,-1)}(Kum)\rightarrow\mathbf{M}_{(0,2\Theta,-2)}^{0}\]
 (defined in Proposition \ref{2:1}) commutes with the $\mathcal{J}[2]$-actions
on $U_{(0,d^{*}H,-1)}(Kum)$ and $\mathbf{M}_{(0,2\Theta,-2)}^{0}$
 described in 2) and 3) of this remark.\\
Moreover, since by 2) and 4) the open subvarieties
$V^{s}\subset\mathbf{M}_{(0,2\Theta,-2)}^{0}$ and
$U^{s}(Kum):=\tau^{-1}(V^{s})$ are $\mathcal{J}[2]$ invariant, the
restriction of $  \tau  $
\[  \tau  ^{s}:U^{s}(Kum)\rightarrow V^{s}\]
defined in \ref{defuskum}  commutes with the induced actions on
$U^{s}(Kum)$ and $V^{s}$.}
\end{enumerate}\end{rem}
 We are now ready to describe the
desired extension of $  \tau  ^{s}$.
\begin{prop}\label{completamento2a1}\begin{enumerate}\item{
Let $j:V^{s}\rightarrow
\widetilde{V}:=\widetilde{\mathbf{M}}_{(0,2\Theta,-2)}^{0}
\setminus\widetilde{\Omega}_{(0,2\Theta,-2)}$ be the natural
inclusion (recall from the definition \ref{def V e Vs} that
$V^{s}$ is a subset of the stable locus of
$\mathbf{M}_{(0,2\Theta,-2)}^{0}$ ), then the
$\mathcal{J}[2]$-action on $V^{s}$ extends to $\widetilde{V}$.}
\item{There exist a smooth algebraic variety $\widetilde{U}(Kum)$,
an open embedding
\[i:U^{s}(Kum)\rightarrow\widetilde{U}(Kum),\] an action of $\mathcal{J}[2]$
 on $\widetilde{U}(Kum)$  extending the given $\mathcal{J}[2]$-action on
  $U^{s}(Kum)$ and  a proper map $\widetilde{  \tau
}:\widetilde{U}(Kum)\rightarrow\widetilde{V}$ making commutative
the following diagram:
\[
\CD
 U^{s}(Kum) @>i>> \widetilde{U}(Kum) \\
@V  \tau  ^{s}VV @V \widetilde{  \tau  }VV \\
V^{s} @>j>>  \widetilde{V}.\\
\endCD
\]
In particular $\widetilde{  \tau  }$ commutes with the
$\mathcal{J}[2]$ actions on $\widetilde{U}(Kum)$ and
$\widetilde{V}$ .}
\end{enumerate}
\end{prop}
\begin{proof}(1): Since by  O'Grady's construction
$\widetilde{V}$ is simply the blow up of
$\mathbf{M}_{(0,2\Theta,-2)}^{0}
\setminus\Omega_{(0,2\Theta,-2)}^{0}$ along
$\Sigma_{(0,2\Theta,-2)}^{0}\setminus\Omega_{(0,2\Theta,-2)}^{0}$
(see Remark \ref{costruzione og rem}), item 1) follows from the
$\mathcal{J}[2]$ invariance of $\Sigma_{(0,2\Theta,-2)}^{0}$ and
$\Omega_{(0,2\Theta,-2)}^{0}$ proved in 2) of remark \ref{rem
j[2]}.\\
(2): By (1) there is a $\mathcal{J}[2]$ equivariant diagram
\begin{equation}\label{primo}
\CD
 U^{s}(Kum) @>id>> U^{s}(Kum) \\
@V  \tau  ^{s}VV @V  \tau  ^{1}VV \\
V^{s} @>j>>  \widetilde{V}\\
\endCD
\end{equation}
where obviously the horizontal maps are open embeddings and $ \tau
^{1}:=j\circ   \tau  ^{s}$.\\ To replace $  \tau  ^{1}$ with a
proper map let $\Gamma\subset U^{s}(Kum)\times\widetilde{V}$ be
the graph of $  \tau  ^{1}$, consider the inclusion
$U^{s}(Kum)\subset\mathbf{M}_{(0,d^{*}H,-1)}(Kum)$ and let
$\overline{\Gamma(Kum)}\subset\mathbf{M}_{(0,d^{*}H,-1)}(Kum)\times\widetilde{V}$
be the closure of $\Gamma(Kum)$. Let $  \tau
^{2}:\overline{\Gamma(Kum)}\rightarrow \widetilde{V}$ be the
restriction of the projection
$p_{2}:\mathbf{M}_{(0,d^{*}H,-1)}(Kum)\times\widetilde{V}\rightarrow\widetilde{V}$
and $i_{\overline{\Gamma}}:U^{s}(Kum)\rightarrow
\overline{\Gamma(Kum)}$ be the natural open embedding: then,
clearly, the following diagram is commutative:
\begin{equation}
\CD
 U^{s}(Kum) @> i_{\overline{\Gamma}}>>  \overline{\Gamma(Kum)}\\
@V  \tau  ^{1}VV @V  \tau  ^{2}VV \\
\widetilde{V} @>id>>  \widetilde{V}\\
\endCD
\end{equation}
 $  \tau  ^{2}$ is proper since it is the restriction to a
closed subvariety of the projection $p_{2}$ which is proper since
$\mathbf{M}_{(0,d^{*}H,-1)}(Kum)$ is projective.\\ Moreover the
$\mathcal{J}[2]$ action on $U^{s}(Kum)$ extends via
$i_{\overline{\Gamma}}$ to $\overline{\Gamma(Kum)}$.\\ To prove
this notice that, since $  \tau  ^{1}$ is $\mathcal{J}[2]$
equivariant, the $\mathcal{J}[2]$-action on $U^{s}(Kum)$ can be
identified with the restriction to $\Gamma(Kum)$ of the
$\mathcal{J}[2]$ diagonal action on
$U^{s}(Kum)\times\widetilde{V}$. Furthermore the $\mathcal{J}[2]$
action on $U^{s}(Kum)\times\widetilde{V}$ is in its turn, by (3)
of remark \ref{rem j[2]}, the restriction of a diagonal action on
$\mathbf{M}_{(0,d^{*}H,-1)}(Kum)\times\widetilde{V}$. Finally,
being $\Gamma(Kum)$ invariant with respect to this action, its
closure $\overline{\Gamma(Kum)}$ is invariant too, and the
restriction of the $\mathcal{J}[2]$-action on
$\mathbf{M}_{(0,d^{*}H,-1)}(Kum)\times\widetilde{V}$ to
$\overline{\Gamma(Kum)}$ provides the desired extension.\\
By a general result, see \cite{AW97, BM97}, for any finite group
acting on an algebraic variety there always exists in
characteristic 0 an equivariant resolution of singularities: let
then $r:\widetilde{U}(Kum)\rightarrow \overline{\Gamma(Kum)}$ be
such a resolution for the $\mathcal{J}[2]$ action on
$\overline{\Gamma(Kum)}$. Since $i_{\overline{\Gamma}}$ is an open
embedding of a smooth variety, it lifts to an open embedding
$i_{\widetilde{U}^{s}(Kum)}:U^{s}(Kum)\rightarrow
\widetilde{U}(Kum)$ and setting $\widetilde{  \tau  }:=  \tau
^{2}\circ r$ we get a $\mathcal{J}[2]$ equivariant commutative
diagram
\begin{equation}
\CD
U^{s}(Kum) @>i_{\widetilde{U}^{s}(Kum)}>> \widetilde{U}(Kum)\\
@V  \tau  ^{1}VV @V\widetilde{  \tau  }VV \\
\widetilde{V} @>id>>  \widetilde{V}\\
\endCD
\end{equation}
 And joining this diagram with
(\ref{primo}) we complete the proof  Proposition
\ref{completamento2a1}.
\end{proof}
\begin{rem} Since by Corollary \ref{etale 2a1} $  \tau  ^{s}$ is proper, we have
$  \tau  ^{s}(\widetilde{U}(Kum)\setminus U^{s}(Kum))\subset
(\widetilde{V}\setminus V^{s})$ , in  particular the image of the
contracted locus does not intersect $V^{s}$.
\end{rem}
Finally we get the divisibility of
$c_{1}(\widetilde{\Sigma}^{0}_{(0,2\Theta,-2)})$.

 \noindent{\em Proof of Proposition \ref{divis
02-2}.}  Let $C\subset\widetilde{U}(Kum)$ be the subvariety
contracted by $\widetilde{  \tau  }$, namely
\[C:=\{x\in\widetilde{V}:
 \;dim(\widetilde{  \tau  }^{-1}(\widetilde{  \tau  }(x)))>0\}\] set
$\widetilde{V}^{0}:=\widetilde{V}\setminus\widetilde{  \tau  }(D)$
and $\widetilde{U}^{0}(Kum):=\widetilde{  \tau
}^{-1}(\widetilde{V}^{0})$, then $\widetilde{V}^{0}$ and
$\widetilde{U}^{0}(Kum)$  are obviously $\mathcal{J}[2]$ invariant
and the restriction
\[\widetilde{  \tau  }^{0}:\widetilde{U}^{0}(Kum)\rightarrow \widetilde{V}^{0}\] is
$\mathcal{J}[2]$ equivariant. Furthermore, by construction,
$\widetilde{  \tau  }^{0}$ is a  proper map with finite fibers ,
hence it is finite, since $\widetilde{U}^{0}(Kum)$ and
$\widetilde{V}^{0}$ are smooth we deduce that $\widetilde{  \tau
}^{0}$ is flat.\\The theory of
 double coverings  can be  applied to deduce that the class
in $Pic(\widetilde{V}^{0})$ of the branch locus of $\widetilde{
\tau }^{0}$ has a square root.
  Since $codim(\widetilde{\mathbf{M}}_{(0,2\Theta,-2)}^{0}\setminus
\widetilde{V}^{0},\widetilde{\mathbf{M}}_{(0,2\Theta,-2)}^{0})>1$,
 the open embedding $\widetilde{V}^{0}\subset
\widetilde{\mathbf{M}}_{(0,2\Theta,-2)}^{0} $ induces an
isomorphism
$Pic(\widetilde{\mathbf{M}}_{(0,2\Theta,-2)}^{0})\simeq Pic
(\widetilde{V}^{0})$ which will give the existence  in
$Pic(\widetilde{\mathbf{M}}_{(0,2\Theta,-2)}^{0})$ of a square
root of the line bundle associated to the closure in
$\widetilde{\mathbf{M}}_{(0,2\Theta,-2)}^{0}$
 of the branch locus of $\widetilde{  \tau  }^{0}$.\\
By the commutativity of diagram in (2) of Proposition
\ref{completamento2a1}, since $\tau^{s}$ is everywhere [2:1] by
Corollary \ref{etale 2a1}, the ramification locus $R$ of
$\widetilde{  \tau  }^{0}$ and obviously its closure
$\overline{R}$ in $\widetilde{\mathbf{M}}^{0}_{(0,2\Theta,-2)}$
are included in
$\widetilde{\Sigma}^{0}_{(0,2\Theta,-2)}\bigcup_{x\in\mathcal{J}[2]}
\Delta_{x}$. Furthermore by the $\mathcal{J}[2]$ equivariance of
$\widetilde{
 \tau  }^{0}$ and since $\mathcal{J}[2]$ acts transitively on the
$\Delta_{x}$, as preannounced, if there is an $x\in \mathcal{J}[2]$ such that
$\Delta_{x}\subset\overline{R}$ then $\bigcup_{x\in\mathcal{J}[2]}
\Delta_{x}\subset\overline{R}$.\\
 Therefore there are only four
possible cases:
\begin{enumerate}
\item{$\overline{R}=\emptyset$,}
\item{$\overline{R}=\bigcup_{x\in\mathcal{J}[2]}\Delta_{x}$,}
\item{$\overline{R}=\widetilde{\Sigma}^{0}_{(0,2\Theta,-2)}\bigcup_{x\in\mathcal{J}[2]}\Delta_{x}$,}
\item{$\overline{R}=\widetilde{\Sigma}^{0}_{(0,2\Theta,-2)}$.}
\end{enumerate}
The first 2 cases are not possible.\\
Suppose by absurd the contrary: then, denoting by $\overline{R}$
also the associated reduced divisor, by Proposition \ref{mult} we
would have the following equality of divisors
\[\overline{R}=\widetilde{\Phi}^{*} D
\] where
$D=0$ if $\overline{R}=\emptyset$, or $D=\sum_{x\in\mathcal{J}[2]}
H_{x}$ if $\overline{R}=\bigcup\Delta_{x}$.\\
 In both of these
cases there would exist a divisor $D^{1}\in Div(|2\Theta|^{\vee})$
such that $[2D^{1}]=[D]$ in $Pic(|2\Theta|^{\vee})$ and since
$Pic(\widetilde{V}^{0})=Pic(\widetilde{\mathbf{M}}_{(0,2\Theta,-2)}^{0})$
is free, the restriction of the class
$[\widetilde{\Phi}^{*}D^{1}]$ to $\widetilde{V}^{0}$ would be the
unique square root of the class $[R]\in Pic(\widetilde{V}^{0})$.
Moreover since the codimension of the complement of
$\widetilde{V}^{0}$ in
$\widetilde{\mathbf{M}}_{(0,2\Theta,-2)}^{0}$ is bigger than $1$,
there would exist a unique (up to scalars) regular section
$\sigma$ of $\mathcal{O}_{\widetilde{V}^{0}}(R)$ vanishing with
multiplicity 1 on each component of $R$ and obviously it would be
the restriction to $\widetilde{V}^{0}$ of the pull back via
$\widetilde{\Phi}$ of a
section $s$ of $\mathcal{O}_{|2\Theta|^{\vee}}(D)$.\\
Therefore the double covering $\widetilde{  \tau  }^{0}$ would be
the one defined by means of $[\widetilde{\Phi}^{*}D]$ and
$\widetilde{\Phi}^{*}s$, in particular its restriction to a fiber
of $\widetilde{\Phi}$ over a $p$ not in $D$ would yield a double
covering defined by means of the trivial square root of the
trivial line bundle and of a constant section, namely we would get
the trivial double covering: this is absurd because for $C\in
|2\Theta|$ smooth, $(\tau^{s})^{-1}(\Phi^{-1}(C))$ is identified
with $M_{C_{K}}\simeq Pic^{1}(C_{K})$, where $C_{K}=f(C)$ is a
smooth genus 3 curve
(see the notation of the proof of Proposition \ref{2:1}).\\
Since by Proposition \ref{mult}
$\mathcal{O}(\sum_{\mathcal{J}[2]}\Delta_{x})=\widetilde{\Phi}^{*}\mathcal{O}(16)$,
both
 $[\overline{R}]$ and
$[\sum_{\mathcal{J}[2]}\Delta_{x}]$ have a square root, therefore
in the last 2 cases $[\widetilde{\Sigma}^{0}_{(0,2\Theta,-2)}]$
has a square root too. \qed

\subsection{The integral basis} We can now present a basis for
$H^{2}(\widetilde{M},\mathbb{Z})$.
\begin{thm}\label{base intera}
Let $\{\alpha_{i}\}_{i=1}^{6}$ be a basis for
$H^{2}(\mathcal{J},\mathbb{Z})$, then  $\{\mu(\alpha_{i}),
c_{1}(\widetilde{B}), A\}$ is a basis for
$H^{2}(\widetilde{M},\mathbb{Z})$.
\end{thm}
\begin{proof}
By Poincar\'e duality, since we already know
$rk(H^{2}(\widetilde{\mathcal{M}},\mathbb{Z}))=8$, it is enough to
show there are 8 elements in
$H_{2}(\widetilde{\mathcal{M}},\mathbb{Z})$ such that the
determinant of the evaluation matrix of these two 8-tuples is 1.
By Proposition \ref{sat} and Poincar\'e duality  we can find
$\{\alpha_{i}^{*}\}_{i=1}^{6}$ such that
$det(<\widetilde{\mu}({\alpha_{i}}),\alpha_{j}^{*}>)=1$. Let
$\gamma$ and $\delta$ be as in the proof of Proposition (7.3.3) of
in \cite{OG03}, then, as shown there we have
\[
\left(
\begin{array}{cc}
<A,\delta> & <A,\gamma> \\
<c_{1}(\widetilde{B}),\delta>  & <c_{1}(\widetilde{B}),\gamma>\\
\end{array}
\right)=
\left(
\begin{array}{cc}
<\frac{c_{1}(\widetilde{\Sigma})}{2},\delta> & <\frac{c_{1}(\widetilde{\Sigma})}{2},\gamma> \\
<c_{1}(\widetilde{B}),\delta>  & <c_{1}(\widetilde{B}),\gamma>\\
\end{array}
\right) = \left(
\begin{array}{cc}
-1 & 1\\
1 & -2\\
\end{array}
\right).
\]
The second equality is proved in Proposition (7.3.3) of in
\cite{OG03}. Furthermore, since $\delta$ and $\gamma$ are
contracted by $\varphi\circ\widetilde{\pi}$, we have
$<\widetilde{\mu}(\alpha_{i}),\delta>=<\widetilde{\mu}(\alpha_{i}),\gamma>=0$
for each i: so the determinant of the whole intersection matrix is given by
\[det(<\widetilde{\mu}(\alpha_{i}),\alpha_{j}^{*}>)\cdot
det\left(\begin{array}{cc}
<A,\delta> & <A,\gamma> \\
<c_{1}(\widetilde{B}),\delta>  & <c_{1}(\widetilde{B}),\gamma>\\
\end{array}
\right)=1\]

\end{proof}
\subsection{Explicit computation}

We now compute the Beauville form of $\widetilde{\mathcal{M}}$ in
terms of the given basis  of
$H^{2}(\widetilde{\mathcal{M}},\mathbb{Z})$.\\
The final result is the following

\begin{thm}\label{mainthm}
Set
$\Lambda:=\mathbb{Z}A\oplus\mathbb{Z}c_{1}(\widetilde{B})\subset
H^{2}(\widetilde{\mathcal{M}},\mathbb{Z})$. There is a direct sum
decomposition
\[H^{2}(\widetilde{\mathcal{M}},\mathbb{Z})=
\widetilde{\mu}(H^{2}(\mathcal{J},\mathbb{Z}))\oplus_{\bot}
\Lambda\] orthogonal with respect to
$B_{\widetilde{\mathcal{M}}}$.\\
The  map
$\widetilde{\mu}:(H^{2}(\mathcal{J},\mathbb{Z}),(,)_{\mathcal{J}})\longrightarrow
(H^{2}(\widetilde{\mathcal{M}},\mathbb{Z}),B_{\widetilde{\mathcal{M}}})$
is an isometric embedding.\\
Furthermore the matrix of the Beauville's form on $\Lambda$  is
given the following formula:
\[
\begin{array}{|c|c|c|}
\hline
& A & c_{1}(\widetilde{B}) \\
\hline
A  & -2 & 2 \\
\hline
c_{1}(\widetilde{B})  & 2 & -4 .\\
\hline
\end{array}
\]
Finally the Fujiki constant of $\widetilde{\mathcal{M}}$ is
$c_{\widetilde{\mathcal{M}}}=60.$
\end{thm}
We will prove this theorem at the end of this section, after
having computed the necessary intersection numbers on
$\widetilde{\mathcal{M}}$.\\ Few intersection numbers are actually
needed to completely determine the Beauville's form on
$\widetilde{\mathcal{M}}$: most of them, those involving
$\widetilde{B}$ and $\widetilde{\Sigma}$, can be computed in terms
of the known geometry of explicit subvarieties ($\widetilde{B}$,
$\widetilde{\Sigma}$ and $\widetilde{B}\cap\widetilde{\Sigma}$) of
$\widetilde{M}$.\\
 Nevertheless to obtain the right normalization
we will need to compute also at least a non zero sextuple
self-intersection of a class in $H^{2}(\widetilde{M})$. More
precisely we want to prove the following proposition.
\begin{prop}\label{muomega6}
\begin{equation}\int_{\widetilde{M}}\widetilde{\mu}(\omega+\overline{\omega})^{6}
=60\left(\int_{\mathcal{J}}(\omega+\overline{\omega})^{2}\right)^{3}.
\end{equation}
\end{prop}This self intersection can be  related  to a known intersection
number on $ Hilb^{3}(Kum)$. The existence of a relation is
provided by the following proposition.
\begin{prop}\label{birHM}
There exists a generically injective rational map
 \[\beta : Hilb^{3}(Kum)
 \dashrightarrow\mathbf{M}_{(0,d^{*}H,-1)}.
 \]
\end{prop}
\begin{proof}
$\beta$ is given as follows. Given 3 general points on $Kum$ there
exists a unique curve $C$ in $|d^{*}H|$ passing through them.
Considering the push forward on $Kum$ of the ideal of the 3 points
on $C$ and tensoring it with $\mathcal{O}(d^{*}H)$, we get a sheaf
in $\mathbf{M}_{(0,d^{*}H,-1)}$. $\beta$ is easily seen to be
generically injective (see \cite{De99,Be99}).
\end{proof}
Therefore we can relate $ Hilb^{3}(Kum)$ with
$\widetilde{\mathcal{M}}$ by means of the rational generically
[2:1] map
\begin{equation}\label{comp}h:=\widetilde{\pi}^{-1}\circ (t^{0})^{-1}\circ (fm^{0})^{-1}\circ   \tau  \circ
\beta:
Hilb^{3}(Kum)\dashrightarrow\widetilde{\mathbf{M}}_{(0,2\Theta,-2)}^{0},
\end{equation}
where $fm^{0}:\mathbf{M}_{(2,2\Theta,0)}^{0}\rightarrow
\mathbf{M}_{(0,2\Theta,-2)}^{0}$ and
$t^{0}:\mathbf{M}_{(2,0,-2)}^{0}\rightarrow
\mathbf{M}_{(2,2\Theta,0)}^{0}$ are  the restrictions of the map
$fm$ and $t$ defined in the proof of Proposition \ref{v02-2}.\\
In order to compute $(\widetilde{\mu}(\omega +
\overline{\omega}))^{6}$ we are going to identify
$h^{*}\circ\widetilde{\mu}(\omega)$. A great simplification in
this identification is provided by the following.
\begin{rem}
The class $\widetilde{\mu}(\omega)\in
H^{2}(\mathcal{M},\mathbb{Z})$ has a holomorphic representative:
this is a particular instance of the last remark on page 504 of
\cite{OG03} asserting that $\widetilde{\mu}$ is a morphism of
Hodge structures.
\end{rem}

Since $h$ is a rational map between smooth varieties it is defined
in codimension 1, by the previous remark
$h^{*}\circ\widetilde{\mu}(\omega)$ is, where defined, an
holomorphic form, therefore by Hartog's theorem it extends to the
whole $ Hilb^{3}( Kum)$. We want to compare this extension that we
will still call $h^{*}\circ\widetilde{\mu}(\omega)$ with another
holomorphic form defined by means of $\omega$ on $ Hilb^{3}
(Kum)$. Precisely, letting $\omega_{K}$ be the unique form on $
Kum$ such that its pull back to $\mathcal{J}$ via the rational
[2:1] map extends to $\omega$, we want to compare
$h^{*}\circ\widetilde{\mu}(\omega)$ with the holomorphic 2-form
$\vartheta({\omega_{K}})$ determined by requiring that its
pull-back to $Kum^{3}$ (via the natural rational map) is equal to
$\sum_{i=1}^{3}\pi_{i}^{*}\omega_{K}$ ($\pi_{i}$ being the i-th
projection).
\begin{thm}\label{rel}
\begin{equation}h^{*}\circ\widetilde{\mu}(\omega)=2\vartheta(\omega_{K})\end{equation}
\end{thm}
 We need a lemma.
\begin{lem}\label{lem1} There is a Zarisky open subset
$U(Kum)\subset Hilb^{3} (Kum)$ such that the following equality of
cohomology classes on $U(Kum)$ holds
\[(h^{*}\circ\widetilde{\mu}(\omega))_{|U(Kum)}=2\vartheta(\omega_{K})_{|U(Kum)}.\]
\end{lem}
\begin{proof}The first step in proving Lemma \ref{lem1}
consists in showing that there exists an open subscheme
$U(Kum)\subset  Hilb^{3}(Kum)$ such that $h$ is defined on
$U(Kum)$, $h(U(Kum))$ is included in the stable locus of
$\mathbf{M}^{0}_{(2,0,-2)}$ and moreover the restriction
\[h_{U}:U(Kum)\rightarrow h(U(Kum))\subset \widetilde{\mathbf{M}}^{0}_{(2,0,-2)}\]
is identified to a modular map induced by a flat family
$\mathcal{F}$ on $\mathcal{J}\times U(Kum)$.\\
 The open subscheme
$U(Kum)$ is simply the locus  $U(Kum)\subset  Hilb^{3}(Kum)$
parametrizing sheaves having a good behavior with respect to the
functors used in the definition of the map $h$. Precisely $U(Kum)$
parametrizes sheaves of ideals $I_{Z}$, where
$Z=\{z_{1},z_{2},z_{3}\}$ consists of 3 distinct points such that:
\begin{enumerate}
\item{Exists a unique $C_{K}\in |d^{*}H|$ such that $z_{i}\in
C_{K}$ and moreover $C_{K}$ is smooth.} \item{Letting
$j:C_{K}\rightarrow Kum$ be the closed embedding  and setting
$p_{i}:=j^{-1}(z_{i})$, the sheaf $f^{*}\circ d_{*}
(j_{*}(\mathcal{O}_{C_{K}}(-\sum p_{i}))\otimes
\mathcal{O}(d^{*}H))$ satisfies W.I.T. (index 1) and its
Fourier-Mukai transform is stable.}
\end{enumerate}
The family $\mathcal{F}$ is defined as follows.\\ Letting
$\mathcal{I}$ be a tautological family of ideals on $Kum\times U(Kum)$ and
letting $p_{K}$ and $q_{K}$ be the projections of $Kum\times
U(Kum)$, by 1) $q_{K*}\mathcal{H}om(p_{K}^{*}\mathcal{O}_{Kum}(-d^{*}H),\mathcal{I})$
is a line bundle
$L$ on $U(Kum)$: the line bundle $q_{K*}\mathcal{H}om(p_{K}^{*}\mathcal{O}_{Kum}(-d^{*}H)
\otimes
q_{K}^{*}L),\mathcal{I})$ has then a nowhere vanishing section
that induces  an exact sequence:
\begin{equation}\label{Definizione G} 0\rightarrow
p_{K}^{*}(\mathcal{O}_{Kum}(-d^{*}H))\otimes q_{K}^{*} L \rightarrow
\mathcal{I}\rightarrow \mathcal{G}\rightarrow 0
\end{equation}
defining a family $\mathcal{G}$ such that $\mathcal{G}_{|Kum\times
Z}\simeq j_{*}(\mathcal{O}_{C_{K}}(-\sum p_{i}))$.\\ Setting
$\mathcal{T}:= (b \times id)_{*}\circ (q \times id)^{*}:
Coh(Kum\times U(Kum))\rightarrow Coh (\mathcal{J}\times U(Kum))$
we can define a new family on $\mathcal{J}\times U(Kum)$
\begin{equation}\label{definizione H}
\mathcal{H}:=\mathcal{T}(\mathcal{G}\otimes
p_{K}^{*}(\mathcal{O}(d^{*}H))).\end{equation} Finally, letting
$p_{J}$ and $q_{J}$ be the projections of $\mathcal{J}\times
U(Kum)$ and denoting $\mathcal{FM}_{U}$  the Fourier-Mukai
transform for families of sheaves on $\mathcal{J}$ parametrized by
$U(Kum)$, using 2) we can define the family $\mathcal{F}$
\begin{equation}\label{Definizione F}
\mathcal{F}:=
p_{J}^{*}(\mathcal{O}(-\Theta))\otimes\mathcal{FM}_{U}(\mathcal{H}).
\end{equation}
and by constrction, the modular map $h_{U}$ induced by
$\mathcal{F}$ is identified on $U(Kum)$ with  $h$. We can now use
the property (\ref{propdonalson}) to compare
$\vartheta(\omega)_{|U(Kum)}$ and $
h^{*}\circ\widetilde{\mu}(\omega)_{|U(Kum)}$. Using
(\ref{propdonalson}) and the orthogonality of $c_{1}(\Theta)$ and
$\omega$:
\begin{equation}\label{primoconto}
h^{*}\circ\widetilde{\mu}(\omega)_{|U(Kum)}= h_{U}^{*}\circ
\varphi^{*}\circ \mu(\omega)=q_{J*}(p_{J}^{*}(\omega)\cup
c_{2}(\mathcal{F}))\end{equation}
\[ =
 q_{J*}(p_{J}^{*}(\omega)\cup c_{2}(\mathcal{FM}_{U}(\mathcal{H})))
=-[q_{J*}(p_{J}^{*}(\omega)\cup
ch(\mathcal{FM}_{U}(\mathcal{H})))]_{2}\] (the last equality holds
since $\mathcal{FM}_{U}(\mathcal{H})$ parametrizes sheaves with
fixed determinant). Decomposing
$ch(\mathcal{FM}_{U}(\mathcal{H}))$ by means of the K\"unneth
decomposition of $H^{\bullet}(\mathcal{J}\times U(Kum))$ and the
Hodge decomposition of $H^{\bullet}(\mathcal{J})$ it is easily
seen that the unique component giving non zero contribution in the
last term of (\ref{primoconto}) is the one in
$H^{0,2}(\mathcal{J})\otimes H^{2}(U(Kum))$. On the other hand
using Grothendieck-Riemann-Roch
\begin{equation}\label{FMGRR}
ch(\mathcal{FM}_{U}(\mathcal{H}))=FM_{U}(ch(\mathcal{H}))
\end{equation}
where $FM_{U}$ is the automorphism of
$H^{\bullet}(\mathcal{J}\times U(Kum))$ induced by means of the
K\"unneth decomposition of $H^{\bullet}(\mathcal{J}\times U(Kum))$
and the automorphism $FM$ induced by $\mathcal{FM}$ on
$H^{\bullet}(\mathcal{J})$ (see \cite{Mu87}). Since
$FM(H^{2}(\mathcal{J}))=H^{2}(\mathcal{J})$ and $FM$, being
defined by means of an algebraic cycle, respects the Hodge
decomposition we get
$FM(H^{0,2}(\mathcal{J}))=H^{0,2}(\mathcal{J})$ and since $FM\circ
FM=-id^{*}$, we also see that $FM$ is the identity on
$H^{0,2}(\mathcal{J})$.\\
It follows that the component in $H^{0,2}(\mathcal{J})\otimes
H^{2}(U(Kum))$ of $ch(\mathcal{FM}_{U}(\mathcal{H}))$ is the same
as the one of $ch(\mathcal{H})$. Therefore
\begin{equation}\label{conto con FM}
[q_{J*}(p_{J}^{*}(\omega)\cup
ch(\mathcal{FM}_{U}(\mathcal{H})))]_{2}=[q_{J*}(p_{J}^{*}(\omega)\cup
ch(\mathcal{H}))]_{2}
\end{equation}
We want now to compare $q_{J*}(p_{J}^{*}(\omega)\cup
ch_{2}(\mathcal{H}))$ and $q_{K*}(p_{K}^{*}(\omega_{K})\cup
ch_{2}(\mathcal{G}\otimes p_{K}^{*}(\mathcal{O}(d^{*}H))))$
recalling by the definition (\ref{definizione H}) that
\[\mathcal{T}(\mathcal{G}\otimes p_{K}^{*}(\mathcal{O}(d^{*}H)))=
(b\times id)_{*}\circ (q\times id)^{*}((\mathcal{G}\otimes
p_{K}^{*}(\mathcal{O}(d^{*}H))))=\mathcal{H}.\] Letting the
algebraic cycle $D:=\sum_{i}n_{i}D_{i}$ be a representative  of
$ch_{2}((\mathcal{G}\otimes p_{K}^{*}(\mathcal{O}(-d^{*}H))))$ not
intersecting  $\bigcup_{x\in\mathcal{J}[2]} E_{x}\times U(Kum)$
and setting $B_{i}:=(b\times id)_{*}\circ (q\times
id)^{*}(D_{i})$, it is easily seen that there exist  \'etale
double coverings $g_{i}:B_{i}\rightarrow D_{i}$ such that
$g_{i}^{*}(p_{K}^{*}(\omega_{K})_{|D_{i}})=(p_{J}^{*}(\omega))_{|B_{i}}$:
since  $\sum_{i}B_{i}$ is a representative of
$ch_{2}(\mathcal{H})$ we get
\begin{equation}\label{conto del 2}
q_{J*}(p_{J}^{*}(\omega)\cup
ch_{2}(\mathcal{H}))=\sum_{i=1}^{n}n_{i}q_{J*}((p_{J}^{*}(\omega))_{|B_{i}})=
\end{equation}\[
2\sum_{i=1}^{n}n_{i}q_{K*}((p_{K}^{*}(\omega_{K}))_{|D_{i}})=2q_{K*}(p_{K}^{*}(\omega_{K})\cup
ch_{2}(\mathcal{G}\otimes p_{K}^{*}(\mathcal{O}(d^{*}H)))) .\]
Finally using the orthogonality of $\omega_{K}$ and $d^{*}(H)$ and
applying Whitney's formula to the sequence (\ref{Definizione G}) a
straightforward computation yields
\begin{equation}\label{ultimoconto}2q_{K*}(p_{K}^{*}(\omega_{K})\cup
ch_{2}(\mathcal{G}\otimes
p_{K}^{*}(\mathcal{O}(d^{*}H))))=2q_{K*}(p_{K}^{*}(\omega_{K})\cup
ch_{2}(\mathcal{G}))=
\end{equation}
\[2q_{K*}(p_{K}^{*}(\omega_{K})\cup
(ch_{2}(\mathcal{I})-ch_{2}(p_{K}^{*}(\mathcal{O}(-d^{*}H))\otimes
q_{K}^{*} L)))=2q_{K*}(p_{K}^{*}(\omega_{K})\cup
ch_{2}(\mathcal{I}))=\]
\[-2\vartheta(\omega_{K})_{|U(Kum)}.\]
The last equality follows from the known description of the 2-cohomology of $Hilb^{3}(Kum)$ (see \cite{OG01}).
 Joining (\ref{ultimoconto})  with the equations
(\ref{primoconto}), (\ref{conto con FM}), (\ref{conto del 2}) and
 we get the Lemma.
\end{proof}

\noindent{\em Proof of Theorem \ref{rel}}. In general, if  $X$
is a projective variety and $U\subset X$ is a Zarisky open
subset, then the natural map  from $H^{0}(\Omega^{2}_{X})$ to
$H^{2}(U,\mathbb{C})$ is injective.
Since by Lemma \ref{lem1} $
h^{*}\circ\widetilde{\mu}(\omega)_{|U}=2\vartheta(\omega_{K})_{|U}
$, it follows that $
h^{*}\circ\widetilde{\mu}(\omega)=2\vartheta(\omega_{K}) $.
\qed\medskip\\
 As an immediate consequence we can compute
$\int_{\widetilde{\mathcal{M}}}\widetilde{\mu}(\omega+\overline{\omega})^{6}$
and, hence,  prove Proposition \ref{muomega6}.\\ \noindent{\em
Proof of Proposition \ref{muomega6}}. By the previous theorem
\[\int_{\widetilde{\mathcal{M}}}\widetilde{\mu}(\omega+\overline{\omega})^{6}=
\frac{1}{2}2^{6}\int_{ Hilb^{3}
(Kum)}\vartheta(\omega+\overline{\omega})^{6}
\]
and by the description of the Beauville form of $Hilb^{3} (Kum)$
(see Theorem (4.2.2) of \cite{OG01})
\[\frac{1}{2}2^{6}\int_{ Hilb^{3} (Kum)}\vartheta(\omega+\overline{\omega})^{6}=
2^{5}\frac{6!}{3!2^{3}}\left(\int_{ Kum}
(\omega_{K}+\overline{\omega}_{K})^{2}\right)^{3}
\]
and since the pull back to $\mathcal{J}$ of $\omega_{K}$ is
$\omega$
\[2^{5}\frac{6!}{3!2^{3}}\left(\int_{ Kum}
(\omega_{K}+\overline{\omega}_{K})^{2}\right)^{3} =
2^{5}\frac{6!}{3!2^{3}}\frac{1}{2^{3}}\left(\int_{\mathcal{J}}(\omega+\overline{\omega})^{2}\right)^{3}=
60\left(\int_{\mathcal{J}}(\omega+\overline{\omega})^{2}\right)^{3}.\]
\qed\medskip

We now pass to compute other intersection numbers of
$\widetilde{\mathcal{M}}$ involving $c_{1}(\widetilde{\Sigma})$
and $c_{1}(\widetilde{B})$. Our first target is the following
proposition.

\begin{prop}\label{int} Let $\alpha_{i}$ be cohomology classes in
$H^{2}(\mathcal{J},\mathbb{Z})$, then
\begin{equation}
\int_{\widetilde{\mathcal{M}}} c_{1}(\widetilde{B})\wedge
c_{1}(\widetilde{\Sigma})
\wedge\widetilde{\mu}(\alpha_{1})\wedge\widetilde{\mu}(\alpha_{2})
\wedge\widetilde{\mu}(\alpha_{3})\wedge\widetilde{\mu}(\alpha_{4})=\;\end{equation}
\[\:2^{4}[(\alpha_{1},\alpha_{2})\cdot(\alpha_{3},\alpha_{4})+
(\alpha_{1},\alpha_{3})\cdot(\alpha_{2},\alpha_{4})+
(\alpha_{1},\alpha_{4})\cdot(\alpha_{2},\alpha_{3})] .\]
\end{prop}
Before proving this proposition we fix the notation.
\begin{notation}Let $\mathcal{P}$ be the Poincar\'e line bundle on
$\mathcal{J}\times\mathcal{\widehat{J}}$, let
$FM:H^{2}(\mathcal{J},\mathbb{Z})$ be the isometry induced by the
Fourier-Mukai transform (see \cite{Mu87}), let $\Delta$ and
$\overline{\Delta}$ be the diagonal and the anti-diagonal on
$\mathcal{J}\times\mathcal{J}$, let $I_{\Delta}$ and
$I_{\overline{\Delta}}$ be their respective  sheaves of ideals :
we can construct on
$\mathcal{J}\times\mathcal{\widehat{J}}\times\mathcal{J}$ the flat
family of sheaves $p_{1,2}^{*}\mathcal{P}\otimes
p_{1,3}^{*}I_{\Delta}\oplus p_{1,2}^{*}\mathcal{P}^{\vee}\otimes
p_{1,3}^{*}I_{\overline{\Delta}}$ (here and in the following
proposition we denote $p_{I}$ the projections of
$\mathcal{J}\times\mathcal{\widehat{J}}\times\mathcal{J}$, the
multiindex $I$ will indicate the image of the projections ). Via
the identification
$\mathcal{J}\leftrightarrow\mathcal{\widehat{J}}$ we can see it as
a family of sheaves on $\mathcal{J}$ parametrized by
$\mathcal{J}\times\mathcal{J}$.\end{notation} The next proposition
describes the pull backs of the images of the Donaldson's morphism
via the modular map
$f:\mathcal{J}\times\mathcal{J}\rightarrow\mathbf{M}_{(2,0,-2)}^{0}$
associated to this family.
\begin{prop}\label{calcoli}Let $\pi_{i}$ be the  projection
of $\mathcal{J}\times\mathcal{J}$ to the i-th factor, then
\[f^{*}\circ\varphi\circ\mu=-2(\pi_{1}^{*}\circ FM+\pi_{2}^{*}).\]
\end{prop}
\begin{proof}
Using the property  (\ref{propdonalson})  of the Donaldson's
morphism we get
\[f^{*}\circ \varphi^{*}\circ \mu (\alpha)=p_{2,3*}(p_{1}^{*}\alpha
\cup c_{2}(p_{1,2}^{*}\mathcal{P}\otimes
p_{1,3}^{*}I_{\Delta}\oplus p_{1,2}^{*}\mathcal{P}^{\vee}\otimes
p_{1,3}^{*}I_{\overline{\Delta}})).\] Since our family
parametrizes sheaves having trivial determinant bundles, and
moreover we have
$c_{1}(I_{\Delta})=c_{1}(I_{\overline{\Delta}})=0$ and
$ch_{2}(\mathcal{P})=ch_{2}(\mathcal{P}^{*})$ it follows
\begin{equation}\label{5}
f^{*}\circ \varphi^{*}\circ \mu
(\alpha)=-p_{2,3*}(p_{1}^{*}\alpha\cup
(2ch_{2}(p_{1,2}^{*}\mathcal{P})-ch_{2}(p_{1,3}^{*}I_{\Delta})-ch_{2}(p_{1,3}^{*}
I_{\overline{\Delta}}))=\end{equation}
\[-p_{2,3*}(p_{1}^{*}\alpha\cup
(2p_{1,2}^{*}(ch_{2}(\mathcal{P}))))-2\pi_{2}^{*}(\alpha)=-
2\pi_{1}^{*}\circ FM (\alpha)-2\pi_{2}^{*}(\alpha);\] the second
equality is verified because  $\Delta$ and $\overline{\Delta}$ act
on the 2-cohomology as the identity ($-id^{*}$ is the identity on
$H^{2}$).
\end{proof}
This proposition enables us  to prove Proposition \ref{int}.\\
\noindent{\em Proof of Proposition \ref{int}}. The modular map $f$
obviously factors through the quotient by the involution $-id$
\[q:\mathcal{J}\times \mathcal{\widehat{J}}\rightarrow \mathcal{J}\times
\mathcal{\widehat{J}}/-id,\] its isomorphic image in
$\mathbf{M}_{(2,0,-2)}^{0}$
is just $\Sigma_{(2,0,-2)}^{0}$.\\
Recall from (7.3.5) of \cite{OG03} that
$\widetilde{B}\cap\widetilde{\Sigma}$ provide a rational section
$\psi$ of the map with general fiber $\mathbb{P}^{1}$ given by
\[\widetilde{\pi}_{|\widetilde{\Sigma}}:\widetilde{\Sigma}
\rightarrow \Sigma_{(2,0,-2)}^{0},\] the composition $\psi\circ
q:\mathcal{J}\times\mathcal{\widehat{J}}\rightarrow\widetilde{B}\cap
\widetilde{\Sigma}$ is then only a rational map: a resolution
$\widetilde{q}$ of $\psi\circ q$ yields  the following commutative
diagram:
\[\CD
 \widetilde{\mathcal{J}\times \mathcal{\widehat{J}}}
@>\widetilde{q}>> \widetilde{\Sigma}\cap \widetilde{B} @>\widetilde{i}>>
\widetilde{\mathcal{M}} \\
@VrVV @V\widetilde{\pi}_{|\widetilde{\Sigma}}VV
@V\widetilde{\pi}VV \\
\mathcal{J}\times \mathcal{\widehat{J}} @>q>> \mathcal{J}\times
\mathcal{\widehat{J}}/-1 @>i>>
\mathbf{M}_{(2,0,-2)}^{0}\\
\endCD.
\]

Given $\alpha_{1}, \alpha_{2}, \alpha_{3}, \alpha_{4}\in
H^{2}(\mathcal{J},\mathbb{Z})$,
 since  $q$ and $\widetilde{q}$ are generically [2:1],
we obtain:

\[
\int_{\widetilde{\mathcal{M}}}c_{1}(\widetilde{B})\wedge c_{1}(
\widetilde{\Sigma})\wedge\bigwedge_{i=1}^{4}\widetilde{\mu}(\alpha_{i})=
\int_{\widetilde{B}\cap\widetilde{\Sigma}}\bigwedge_{i=1}^{4}\widetilde{i}^{*}\circ
\widetilde{\mu}(\alpha_{i})=
\frac{1}{2}\int_{\widetilde{\mathcal{J}\times
\mathcal{\widehat{J}}}}\bigwedge_{i=1}^{4}
\widetilde{q}^{*}\circ\widetilde{i}^{*}\circ\widetilde{\mu}(\alpha_{i}).\]

By the commutativity of the previous diagram and since $r$ is
birational \[\frac{1}{2}\int_{\widetilde{\mathcal{J}\times
\mathcal{\widehat{J}}}}\bigwedge_{i=1}^{4}
\widetilde{q}^{*}\circ\widetilde{i}^{*}\circ\widetilde{\mu}(\alpha_{i})=
 \frac{1}{2}\int_{\widetilde{\mathcal{J}\times \mathcal{\widehat{J}}}}
\bigwedge_{i=1}^{4}r^{*}\circ q^{*}\circ i^{*}\circ
\varphi^{*}\circ\mu(\alpha_{i})=\]
\[
\frac{1}{2}\int_{\mathcal{J}\times
\mathcal{\widehat{J}}}\bigwedge_{i=1}^{4} q^{*}\circ
i\circ\varphi^{*}\circ\mu(\alpha_{i})=
\frac{1}{2}\int_{\mathcal{J}\times
\mathcal{\widehat{J}}}\bigwedge_{i=1}^{4} f^{*}\circ
\varphi^{*}\circ\mu(\alpha_{i}).\]

Proposition \ref{calcoli} allows us to compute the last term of
this equations by means of the intersection form on
$H^{2}(\mathcal{J},\mathbb{Z})$, we obtain
\[\frac{1}{2}\int_{\mathcal{J}\times
\mathcal{\widehat{J}}}\bigwedge_{i=1}^{4} f^{*}\circ
\varphi^{*}\circ\mu(\alpha_{i})=\frac{1}{2}\cdot2^{4}
\int_{\mathcal{J}\times
\mathcal{\widehat{J}}}\bigwedge_{i=1}^{4}(\pi_{1}^{*}\circ
FM(\alpha_{i})+\pi_{2}^{*}\alpha_{i}).\] Since $FM$ is an isometry
(\cite{Mu87}) by an easy computation this last term equals
\[\frac{1}{2}\cdot2^{4}\cdot2[(\alpha_{1},\alpha_{2})\cdot(\alpha_{3},\alpha_{4})+
(\alpha_{1},\alpha_{3})\cdot(\alpha_{2},\alpha_{4})+
(\alpha_{1},\alpha_{4})\cdot(\alpha_{2},\alpha_{3})]\] as asserted
in the statement of this proposition. \qed\medskip\\
 The formula of Proposition
\ref{int} will give information about the Beauville form
$B_{\widetilde{\mathcal{M}}}$ of $\widetilde{\mathcal{M}}$ when
compared with the Fujiki polarized formula (\ref{formula fujiki
polariizzata}) for the same integral. As a first consequence we
obtain the following orthogonality relations.
\begin{prop}\label{ortogonalita'}
For $\alpha\in H^{2}(\mathcal{J},\mathbb{Z})$:
\[B_{\widetilde{\mathcal{M}}}(c_{1}(\widetilde{B}),\widetilde{\mu}(\alpha))=
B_{\widetilde{\mathcal{M}}}(c_{1}(\widetilde{\Sigma}),\widetilde{\mu}(\alpha))=0.\]
\end{prop}
\begin{proof}
Let $\alpha\in H^{2}(\mathcal{J},\mathbb{Z})$ be a class such that
$(\alpha,\alpha)\neq 0$, then the Proposition \ref{int} applied in
the case $\alpha_{i}=\alpha$ gives
\[\int_{\widetilde{\mathcal{M}}}c_{1}(\widetilde{B})\wedge c_{1}(\widetilde{\Sigma})
\wedge\widetilde{\mu}(\alpha)^{4}= 2^{4}\cdot 3(\alpha,\alpha)^{2}.\] On
the other hand, Fujiki's polarized formula (\ref{formula fujiki
polariizzata}) applied in the same case gives
\[\int_{\widetilde{\mathcal{M}}}c_{1}(\widetilde{B})\wedge c_{1}(\widetilde{\Sigma})\wedge\widetilde{\mu}(\alpha)^{4}=
a_{1}B_{\widetilde{\mathcal{M}}}(c_{1}(\widetilde{\Sigma}),c_{1}(\widetilde{B}))
B_{\widetilde{\mathcal{M}}}(\widetilde{\mu}(\alpha),\widetilde{\mu}(\alpha))^{2}+
\]\[a_{2}B_{\widetilde{\mathcal{M}}}(c_{1}(\widetilde{\Sigma}),\widetilde{\mu}(\alpha))
B_{\widetilde{\mathcal{M}}}(c_{1}(\widetilde{B}),\widetilde{\mu}(\alpha))
B_{\widetilde{\mathcal{M}}}(\widetilde{\mu}(\alpha),\widetilde{\mu}(\alpha)),\]
($a_{1}$ and $a_{2}$ are suitable constants ) hence we deduce,
under our hypothesis,
$B_{\widetilde{\mathcal{M}}}(\widetilde{\mu}(\alpha),\widetilde{\mu}(\alpha))\neq
0$.
 Moreover, since
$\widetilde{\Sigma}$ and $\widetilde{B}$ are contracted in the
Uhlenbeck compactification \[
\int_{\widetilde{\mathcal{M}}}c_{1}(\widetilde{\Sigma})\wedge
\widetilde{\mu}(\alpha)^{5}=\int_{\widetilde{\mathcal{M}}}c_{1}(\widetilde{B})\wedge
\widetilde{\mu}(\alpha)^{5}=0,\] and, applying again the Fujiki's
polarized formula:
\[B_{\widetilde{\mathcal{M}}}(c_{1}(\widetilde{\Sigma}),\widetilde{\mu}(\alpha))
B_{\widetilde{\mathcal{M}}}(\widetilde{\mu}(\alpha),\widetilde{\mu}(\alpha))^{2}=
B_{\widetilde{\mathcal{M}}}(c_{1}(\widetilde{B}),\widetilde{\mu}(\alpha))B_{\widetilde{\mathcal{M}}}(\widetilde{\mu}(\alpha),\widetilde{\mu}(\alpha))^{2}=0.\]
The statement of the proposition holds also for the general
$\alpha$ since the intersection form of $
H^{2}(\mathcal{J},\mathbb{Z})$ is non degenerate.
\end{proof}
It is now possible to determine, up to a scalar factor, the
restriction of $B_{\widetilde{\mathcal{M}}}$ to the image of the
Donaldson's map.
\begin{prop}\label{a meno di costante}$\exists a\in\mathbb{Q}$ such that for any $\alpha\in
H^{2}(\mathcal{J},\mathbb{Z})$
\[B_{\widetilde{\mathcal{M}}}(\widetilde{\mu}(\alpha),\widetilde{\mu}(\alpha))=a(\alpha,\alpha).\]
\end{prop}
\begin{proof} Applying Proposition \ref{int}
we find
\[(\alpha,\alpha)^{2}=
\frac{1}{2^{4}\cdot 3}\int_{\widetilde{\mathcal{M}}}
c_{1}(\widetilde{B})\wedge
c_{1}(\widetilde{\Sigma})\wedge\widetilde{\mu}(\alpha)^{4}\] and
using the polarized Fujiki's formula (\ref{formula fujiki
polariizzata}), thanks to the just proved orthogonality relations:
\[B_{\widetilde{\mathcal{M}}}
(\widetilde{\mu}(\alpha),\widetilde{\mu}(\alpha))^{2}=
\frac{6!}{4!\cdot6\cdot c_{\widetilde{\mathcal{M}}}\cdot B_{\widetilde{\mathcal{M}}}(c_{1}(\widetilde{B}),
c_{1}(\widetilde{\Sigma}))}
\int_{\widetilde{\mathcal{M}}}c_{1}(\widetilde{B})\wedge
c_{1}(\widetilde{\Sigma})\wedge\widetilde{\mu}(\alpha)^{4}\] and
since their squares are proportional, the restriction of
$B_{\widetilde{\mathcal{M}}}$ to the image of the Donaldson's
morphism and the intersection forms on $\mathcal{J}$ are
proportional.
\end{proof}

We need a few other intersection numbers on
$\widetilde{\mathcal{M}}$ to completely determine the Beauville
form on $H^{2}(\widetilde{\mathcal{M}},\mathbb{Z})$: we collect
them in the following proposition.
\begin{prop}\label{ultime formule}

$\forall\omega\in H^{0}(\Omega^{2},\mathcal{J})$ the followings
hold:
\begin{enumerate}\item{$\int_{\widetilde{\mathcal{M}}}c_{1}(\widetilde{B})
\wedge c_{1}(\widetilde{\Sigma})\wedge
\widetilde{\mu}^{4}(\omega+\overline{\omega})= 2^{4}\cdot3
(\omega+\overline{\omega},\omega+\overline{\omega})^{2}$,}
\item{$\int_{\widetilde{\mathcal{M}}}c_{1}(\widetilde{\Sigma})^{2}\wedge
\widetilde{\mu}^{4}(\omega+\overline{\omega})=-2^{5}\cdot3
(\omega+\overline{\omega},\omega+\overline{\omega})^{2}$.}
\item{$\int_{\widetilde{\mathcal{M}}}c_{1}(\widetilde{B})^{2}\wedge
\widetilde{\mu}^{4}(\omega+\overline{\omega})= -2^{4}\cdot3
(\omega+\overline{\omega},\omega+\overline{\omega})^{2}$,}

\end{enumerate}

\end{prop}
\begin{proof}
The first formula is again a  special case of Proposition
\ref{int}. To prove the others consider the more general case of
an holomorphic $\mathbb{P}^{1}$-bundle

\[
\begin{array}{ccc}
\rho & \stackrel{i_{1}}{\hookrightarrow} & F \\
 &  & p\downarrow \\
 &  & U \\
\end{array}
\]
with a smooth base, non necessarily compact.\\ Suppose $\omega$ be
a holomorphic 2-form on $F$ and $Y\subset F$ be an algebraic
smooth closed
 subvariety such that $\rho\cdot Y=d$:
then there exists  on $U$ a unique holomorphic 2-form $\omega_{U}$
whose pull back to $F$ is $\omega$, and obviously
\[\int_{Y}(\omega+\overline{\omega})^{n}=d\int_{U}(\omega_{U}+\overline{\omega}_{U})^{n}.\]
Applying this remark to our case, since $\widetilde{\Sigma}\cap
\widetilde{B}$ is a rational section of the restriction of
$\varphi\circ\widetilde{\pi}$ to $\widetilde{\Sigma}$ (see page 504 of \cite{OG03}) and by
adjunction the normal bundle to $\widetilde{\Sigma}$ has degree -2
on the fibers of such a fibration, we find
\[\int_{\widetilde{\mathcal{M}}}c_{1}(\widetilde{\Sigma})^{2}\wedge
\widetilde{\mu}(\omega+\overline{\omega})^{4}=
-2\int_{\widetilde{\mathcal{M}}}c_{1}(\widetilde{B})\wedge
c_{1}(\widetilde{\Sigma})\wedge
\widetilde{\mu}(\omega+\overline{\omega})^{4}.\]
Analogously
the restriction of $\varphi\circ\widetilde{\pi}$ to $\widetilde{B}$
is still generically a
$\mathbb{P}^{1}$ fibration, $\widetilde{\Sigma}$ intersect its general fiber in 2 points
(see \cite{OG03}) and the degree of the normal bundle of
 $\widetilde{B}$ on the fiber is $-2$: we deduce
\[\int_{\widetilde{\mathcal{M}}}c_{1}(\widetilde{B})^{2}\wedge
\widetilde{\mu}(\omega+\overline{\omega})^{4}=
-\int_{\widetilde{\mathcal{M}}}c_{1}(\widetilde{B})\wedge
c_{1}(\widetilde{\Sigma})\wedge
\widetilde{\mu}(\omega+\overline{\omega})^{4}.\] These two
formulas and the first of the statement obviously imply the
proposition.
\end{proof}
We are now ready to prove Theorem \ref{mainthm}.\\
\noindent{\em Proof of \ref{mainthm}.} The orthogonality relations
have been proved in Proposition \ref{ortogonalita'}. By Fujiki's
polarized formula (see (\ref{formula fujiki polariizzata})) and
using Proposition \ref{ultime formule}

\[c_{\widetilde{\mathcal{M}}}\frac{4!6}{6!}B_{\widetilde{\mathcal{M}}}
(c_{1}(\widetilde{B}),c_{1}(\widetilde{\Sigma}))
B_{\widetilde{\mathcal{M}}}^{2}
(\widetilde{\mu}(\omega+\overline{\omega}),\widetilde{\mu}(\omega+\overline{\omega}))
 =\int_{\widetilde{\mathcal{M}}}c_{1}(\widetilde{B})\wedge c_{1}(\widetilde{\Sigma})\wedge
\widetilde{\mu}(\omega+\overline{\omega})^{4}=\]\[
2^{4}\cdot3(\omega+\overline{\omega},\omega+\overline{\omega})^{2},\]

\[c_{\widetilde{\mathcal{M}}}\frac{4!6}{6!}B_{\widetilde{\mathcal{M}}}
(c_{1}(\widetilde{\Sigma}),c_{1}(\widetilde{\Sigma}))
B_{\widetilde{\mathcal{M}}}^{2}
(\widetilde{\mu}(\omega+\overline{\omega}),\widetilde{\mu}(\omega+\overline{\omega}))
 =\int_{\widetilde{\mathcal{M}}}c_{1}(\widetilde{\Sigma})^{2}\wedge
\widetilde{\mu}(\omega+\overline{\omega})^{4}=\]\[
-2^{5}\cdot3(\omega+\overline{\omega},\omega+\overline{\omega})^{2}.\]

\[c_{\widetilde{\mathcal{M}}}\frac{4!6}{6!}
B_{\widetilde{\mathcal{M}}}(c_{1}(\widetilde{B}),c_{1}(\widetilde{B}))
B_{\widetilde{\mathcal{M}}}^{2}
(\widetilde{\mu}(\omega+\overline{\omega}),\widetilde{\mu}(\omega+\overline{\omega}))
 =\int_{\widetilde{\mathcal{M}}}c_{1}(\widetilde{B})^{2}\wedge
\widetilde{\mu}(\omega+\overline{\omega})^{4}=\]\[
-2^{4}\cdot3(\omega+\overline{\omega},\omega+\overline{\omega})^{2},\]

Reading only the first and the third terms of these equations,
recalling that $c_{1}(\widetilde{\Sigma})=2A$ and using, by
Proposition \ref{a meno di costante},
$B_{\widetilde{\mathcal{M}}}(\widetilde{\mu}(\omega+\overline{\omega}),
\widetilde{\mu}(\omega+\overline{\omega}))=
a(\omega+\overline{\omega},\omega+\overline{\omega})$ we find

\[B_{\widetilde{\mathcal{M}}}(A,c_{1}(\widetilde{B}))=
\frac{2^{4}\cdot3}{2a^{2}\cdot
c_{\widetilde{\mathcal{M}}}\frac{4!6}{6!}}=\frac{2^{3}\cdot3\cdot
5}{a^{2}\cdot c_{\widetilde{\mathcal{M}}}},\]

\[B_{\widetilde{\mathcal{M}}}(A,A)=
\frac{-2^{5}\cdot3}{4a^{2}\cdot
c_{\widetilde{\mathcal{M}}}\frac{4!6}{6!}}=\frac{-2^{3}\cdot3\cdot
5}{a^{2}\cdot c_{\widetilde{\mathcal{M}}}}.\]

\[B_{\widetilde{\mathcal{M}}}(c_{1}(\widetilde{B}),c_{1}(\widetilde{B}))=
\frac{-2^{4}\cdot3}{a^{2}\cdot
c_{\widetilde{\mathcal{M}}}\frac{4!6}{6!}}=\frac{-2^{4}\cdot3\cdot
5}{a^{2}\cdot c_{\widetilde{\mathcal{M}}}},\]

Furthermore by Corollary \ref{muomega6} and Fujiki's formula
(\ref{formula fujiki}) we get
\[
c_{\widetilde{\mathcal{M}}}B_{\widetilde{\mathcal{M}}}^{3}
(\widetilde{\mu}(\omega+\overline{\omega}),\widetilde{\mu}(\omega+\overline{\omega}))=
\int_{\widetilde{\mathcal{M}}}
\widetilde{\mu}^{6}(\omega+\overline{\omega})=
60(\omega+\overline{\omega},\omega+\overline{\omega})^{3}
\]
which implies $c_{\widetilde{\mathcal{M}}}=\frac{60}{a^{3}}$ and
substituting this in the previous formulas :
\[\left(
\begin{array}{cc}
B_{\widetilde{\mathcal{M}}}(A,A) & B_{\widetilde{\mathcal{M}}}(A, c_{1}(\widetilde{B})) \\
B_{\widetilde{\mathcal{M}}}(c_{1}(\widetilde{B}),A)  &
 B_{\widetilde{\mathcal{M}}}(c_{1}(\widetilde{B}),c_{1}(\widetilde{B}))\\
\end{array}
\right) = \left(
\begin{array}{cc}
-2a & 2a\\
 2a & -4a\\
\end{array}
\right).
\] The  rational number $a$ has then to be integer to
respect the integrality of $B_{\widetilde{\mathcal{M}}}$,
moreover, since $a$ divides $B_{\widetilde{\mathcal{M}}}$, the
primitivity of the Beauville form imposes $|a|=1$ and, finally,
since
$c_{\widetilde{\mathcal{M}}}>0$ then $a=1$.\\
 The Fujiky constant is
then $c_{\widetilde{\mathcal{M}}}= 60$ and the bilinear form
$B_{\widetilde{\mathcal{M}}}$ is the one described in the
statement of this theorem. \qed\medskip\\
By an easy change of base we get:
\begin{cor}\label{unimodularita'}
There exsists an isomorphism of lattices:
\[(H^{2}(\widetilde{\mathcal{M}},\mathbb{Z}),
B_{\widetilde{\mathcal{M}}})\simeq (H^{2}(\mathcal{J},\mathbb{Z}),(,)_{\mathcal{J}})\oplus_{\perp}\mathbb{Z}A
\oplus_{\perp}\mathbb{Z}C\] where
 $A\cdot A=C\cdot C=-2$.
\end{cor}


\begin{thebibliography}{100}
\addcontentsline{toc}{chapter}{Bibliography}
\bibitem[AW 97]{AW97}                \textsc{D.Abramovich, J.Wang}, `Equivariant resolution of singularities in characteristic 0',
                                {\em Math. Res. Lett.}{\bf4} (1997), 427--433.

\bibitem[BDL 01]{BDL01}              \textsc{J.Bryan, R.Donagi, N.C.Leung}, `G-Bundles on abelian surfaces, hyperk\"ahler manifolds, and stringy Hodge numbers',
                                {\em Turkish J. Math} {\bf25} (2001), 195--236.

\bibitem[Be 83]{Be83}          \textsc{A.Beauville}, `Vari\'et\'e K\"ahl\'eriennes dont la premi\`ere classe de
                                chern est nulle',
                                {\em J. Diff. Geom.} {\bf18} (1983).

\bibitem[Be 99]{Be99}          \textsc{A.Beauville}, `Counting rational curves on K3 surfaces',
                                {\em Duke Math. J. } {\bf97} (1999).

\bibitem[BM 97]{BM97}          \textsc{E.Bierston, P.Milman}, `Canonical desingularization in characteristic zero by blowing up the maximum strata of a local invariant',
                                {\em Invent. Math.} {\bf128} (1997), 207--302.

\bibitem[Bo 74]{Bo74}          \textsc{F.Bogomolov}, `On the decomposition of K\"ahler manifolds with trivial canonical class',
                                {\em Math. USSR-Sb} {\bf22} (1974), 580--583.

\bibitem[BPV 84]{BPV84}        \textsc{W.Barth, C.Peters, A.Van de Ven}, `Compact complex surfaces',
                                {\em Springer} (1984)

\bibitem[De 99]{De99}          \textsc{O.Debarre}, `On the Euler characteristic of generalized Kummer varieties',
                                {\em Amer.J.Math} {\bf121} (1999), 577--586.

\bibitem[Ei 95]{Ei95}          \textsc{D.Eisenbud}, `Commutative algebra with a view toward algebraic geometry',
                                {\em Springer GTM} {\bf150} (1995).

\bibitem[FM 94]{FM94}          \textsc{R.Friedman, J.Morgan}, `Smooth four-manifolds and complex surfaces',
                                {\em Ergeb. Math. Grenzgeb. 3. Folge 27 Springer} (1994).

\bibitem[Fu 87]{Fu87}          \textsc{A.Fujiki}, `On the de Rham Cohomology Group  of a Compact  K\"ahler Symplectic Manifold',
                                {\em Adv. Studies in Pure Math.} {\bf10}(1987) 105--165.



\bibitem[HL 97]{HL97}          \textsc{D.Huybrechts, M.Lehn}, `The geometry of moduli spaces of scheaves',
                                {\em Aspects of Mathematics} {\bf E 31}(1997).

\bibitem[HM 98]{HM98}          \textsc{J.Harris, I.Morrison}, `Moduli of curve',
                                {\em Graduate Texts in Mathematics} {\bf187} (1998).

\bibitem[Hu 97]{Hu97}          \textsc{D.Huybrechts}, `Birational symplectic manifolds and their deformations',
                                {\em J. Diff. Geom.} {\bf45} (1997), 488--513.


\bibitem[Li 93]{Li93}          \textsc{J.Li}, `Algebraic geometric interpretation of Donaldson's polynomial invariants of algebraic surfaces',
                                {\em J. Diff. Geom.} {\bf37} (1993), 417--466.

\bibitem[LP 93]{LP93}          \textsc{J. Le Potier}, `Syst\`emes coh\'erents et structures de niveau',
                                {\em Ast\'erisque} {\bf214} (1993).

\bibitem[Ma 99]{Ma99}          \textsc{D.Matsushita}, `On fibre space structures of a projective irreducible symplectic manifolds',
                                {\em Topology} {\bf38} (1999), 79--83.

\bibitem[Ma 01]{Ma01}          \textsc{D.Matsushita}, Addendum to `On fibre space structures of a projective irreducible symplectic manifolds',
                                {\em Topology} {\bf40} (2001), 431--432.

\bibitem[Mo 93]{Mo93}          \textsc{ J. Morgan}, `Comparison of the Donaldson polynomial invariants with their algebro-geometric analogues',
                                {\em Topology} {\bf32} (1993), 449--488.

\bibitem[Mu 81]{Mu81}          \textsc{ S.Mukai}, `Duality between D(X) and D($\widehat{X}$) with its applications to Picard sheaves',
                                {\em Nagoya Math. J.} {\bf81} (1981), 153--175.

\bibitem[Mu 84]{Mu84}          \textsc{ S.Mukai}, `Symplectic structure of the moduli space of sheaves on an abelian or K3 surface',
                                {\em Invent. Math} {\bf77} (1984), 101--116.

\bibitem[Mu 87]{Mu87}          \textsc{S.Mukai}, `Fourier functor and its application to the moduli of bundles on an abelian variety',
                                {\em Adv. Stud. Pure Math.} {\bf10} (1987), 515--550.

\bibitem[NR 69]{NR69}          \textsc{M.S.Narasimhan, S.Ramanan}, `Moduli of vector bundles on a compact Riemann surface',
                                {\em Ann. of Math. } {\bf89} (1969), 14--51.

\bibitem[OG 99]{OG99}          \textsc{ K.O'Grady}, `Desingularized moduli spaces of sheaves on a K3',
                                {\em J.Reine Angew. Math.} {\bf512} (1999), 49--117.

\bibitem[OG 01]{OG01}          \textsc{K.O'Grady}, `Compact Hyperk\"ahler Manifolds',  on the web page
                                www.mat.uniroma1.it/people/ogrady.

\bibitem[OG 03]{OG03}          \textsc{ K.O'Grady}, `A new six dimensional irreducible symplectic variety',
                                {\em J. Algebraic Geometry} {\bf12} (2003), 435--505.


\bibitem[Sa 03]{Sa03}          \textsc{J. Sawon}, `Abelian fibred holomorphic symplectic manifolds',
                                {\em Turkish Journal of Mathematics} {\bf27} (2003), 197--230.

\bibitem[Yo 01]{Yo01}          \textsc{K.Yoshioka}, `Moduli spaces of stable sheaves on abelian surfaces',
                                {\em Math.Ann.} {\bf 321} (2001), 817--884.

\bibitem[Yo]{Yo}                \textsc{ K.Yoshioka}, `Irreducibility of moduli spaces of vector bundles on K3 surfaces',
                                {\em math.AG/9907001}.




\end{thebibliography}
\end{document}